\documentclass[11pt,a4paper]{article}

\usepackage{amsxtra}
\usepackage{amssymb}
\usepackage{amsmath}
\usepackage{latexsym}
\usepackage{amsthm}
\usepackage{epsfig}
\usepackage{color}
\usepackage[dvipsnames]{xcolor}
\usepackage{psfrag}
\usepackage{graphicx}

\usepackage{hyperref}

\newtheorem{theorem}{Theorem}[section]

\newtheorem{lemma}[theorem]{Lemma}
\newtheorem{proposition}{Proposition}

\theoremstyle{definition}
\newtheorem{definition}[theorem]{Definition}

\newcommand{\R}{\mathbb{R}}

\newcommand{\Z}{\mathbb{Z}}

\def\cC{{\mathcal C}}

\def\cS{{\mathcal S}}

\newcommand{\Inf}{\displaystyle \inf}

\newcommand{\Sum}{\displaystyle \sum}

\def\ddl{\dot \Delta_l}
\def\ddq{\dot \Delta_q}
\def\tilde{\widetilde}
\def\hat{\widehat}
\newcommand{\D}{\Delta}

\newcommand{\n}{\nabla}
\newcommand{\fd}{\frac{d}{2}}
\newcommand{\fdp}{\frac{d}{p}}
\newcommand{\p}{\partial}

\newcommand{\de}{\delta}

\newcommand{\NN}{\frac{N}{p}}

\newcommand{\N}{\frac{N}{2}}

\newcommand{\e}{\epsilon}
\newcommand{\va}{\varphi}

\newcommand{\T}{\mathbb{T}}

\newtheorem{remarka}{Remark}

\newtheorem{lemme}{Lemma}

\title{Global strong solution for the Korteweg system with quantum pressure in dimension $N\geq 2$}
\author{Boris Haspot  \thanks{Universit\'e Paris Dauphine, PSL Research University, Ceremade, Umr Cnrs 7534, Place du Mar\' echal De Lattre De Tassigny 75775 Paris cedex 16 (France), haspot@ceremade.dauphine.fr }}
\date{}
\begin{document}
\maketitle

\begin{abstract}
This work is devoted to prove the existence of global strong solution in dimension $N\geq 2$ for a 
isothermal model of capillary fluids derived by J.E Dunn and
J.Serrin (1985)  (see \cite{fDS}), which can be used as a phase transition model. We will restrict us to the case of the so called compressible Navier-Stokes system with quantum pressure which corresponds to consider  the capillary coefficient $\kappa(\rho)=\frac{\kappa_1}{\rho}$ with $\kappa_1>0$. In a first part we prove the existence of strong solution in finite time for large initial data with a precise bound by below on the life span $T^*$. This one depends on the norm of the initial data $(\rho_0,v_0)$. The second part consists in proving the existence of global strong solution with particular choice on the capillary coefficient ( where $\kappa_1=\mu^2$) and on the viscosity tensor which corresponds to the viscous shallow water case $-2\mu{\rm div}(\rho Du)$. To do this we derivate different energy estimate on the density and the effective velocity $v$ which ensures that the strong solution can be extended beyond $T^*$. The main difficulty consists in controlling the vacuum or in other words to estimate the $L^\infty$ norm of $\frac{1}{\rho}$. The proof relies mostly on a method introduced by De Giorgi \cite{DG} (see also Ladyzhenskaya et al in \cite{La} for the parabolic case) to obtain regularity results for elliptic equations with discontinuous diffusion coefficients and a suitable bootstrap argument.
\end{abstract}
\tableofcontents
\section{Introduction}
We are concerned with compressible fluids endowed with internal
capillarity. The model we consider  originates from the XIXth
century work by Van der Waals and Korteweg \cite{VW,fK} and was
actually derived in its modern form in the 1980s using the second
gradient theory, see for instance \cite{fDS,fJL,fTN}. The first investigations begin with the Young-Laplace theory which claims that the phases are separated by a hypersurface and that the jump in the pressure across the hypersurface is proportional to the curvature of the hypersurface. The main difficulty consists in describing the location and the movement of the interfaces.\\
Another major problem is to understand whether the interface behaves as a discontinuity in the state space (sharp interface SI) or whether the phase boundary corresponds to a more regular transition (diffuse interface, DI).
The diffuse interface models have the advantage to consider only one set of equations in a single spatial domain (the density takes into account the different phases) which considerably simplifies the mathematical and numerical study (indeed in the case of sharp interfaces, we have to treat a problem with free boundary).\\
Another approach corresponds to determine equilibrium solutions which classically consists in the minimization of the free energy functional.
Unfortunately this minimization problem has an infinity of solutions, and many of them are physically wrong. In order to overcome this difficulty, Van der Waals in the XIX-th century was the first to add a term of capillarity to select the physically correct solutions, modulo the introduction of a diffuse interface. This theory is widely accepted as a thermodynamically consistent model for equilibria. Alternatively, another way to penalize the high density variations consists in applying a zero order but non-local operator to the density gradient (we refer to \cite{9Ro}, \cite{5Ro}, \cite{Rohdehdr}).  We refer for a mathematical analysis on this system to \cite{CHa,CH, CH1, Has1,Has4}.\\
Let us now consider a fluid of density $\rho\geq 0$, velocity field $u\in\R^N$, we are now interested in the following
compressible capillary fluid model, which can be derived from a Cahn-Hilliard  free energy (see the
pioneering work by Dunn and  Serrin in \cite{fDS} and also in
\cite{fA,fC,fGP}).
The conservation of mass and of momentum write:
\begin{equation}
\begin{cases}
\begin{aligned}
&\frac{\p}{\p t}\rho+{\rm div}(\rho u)=0,\\
&\frac{\p}{\p t}(\rho u)+{\rm div}(\rho
u\otimes u)-\rm div(\mu(\rho)\,Du)-\n(\lambda(\rho){\rm div}u)+\n P(\rho)={\rm div}K,
\end{aligned}
\end{cases}
\label{13systeme}
\end{equation}
where the Korteweg tensor reads as following:
\begin{equation}
{\rm div}K
=\n\big(\rho\kappa(\rho)\D\rho+\frac{1}{2}(\kappa(\rho)+\rho\kappa^{'}(\rho))|\n\rho|^{2}\big)
-{\rm div}\big(\kappa(\rho)\n\rho\otimes\n\rho\big).
\label{divK}
\end{equation}
Here $\kappa$ is the capillary coefficient and is a regular function of the form $\kappa(\rho)=\kappa\rho^{\alpha}$ with $\alpha\in\R$. The term
${\rm div}K$  allows to describe the variation of density at the interfaces between two phases, generally a mixture liquid-vapor. The pressure  $P(\rho)=a\rho^{\gamma}$ with $\gamma\geq 1$ is a general $\gamma$ law pressure term, $\mu(\rho)>0$ and $\lambda(\rho)$ are the viscosity coefficient and $D u=\frac{1}{2}(\n u+^{t}\n u)$ is the strain tensor. 
\begin{remarka}
In the sequel we are focusing on the case of shallow-water viscosity coefficients, it means $\mu(\rho)=2\mu\rho$ with $\mu>0$ and $\lambda(\rho)=0$. In addition we will deal only with the case of the quantum compressible Navier-Stokes system  studied in particular in \cite{J} which corresponds to the capillary coefficient $\kappa(\rho)=\frac{\kappa}{\rho}$. When $\mu=0$ the system corresponds to the Euler system with quantum pressure. In \cite{Antonelli}, the authors prove the existence of global weak solution for irrotational initial data, in \cite{C1,C2} we prove with  Audiard the existence of global strong solution for small irrotational initial data.
\end{remarka}
Let us express now the energy of the system (when $P(\rho)=a\rho$ with $a>0$ for simplicity) and when the density is close from a constant state $\bar{\rho}>1$, multiplying the momentum equation by $u$ we have:
with $q(t)=t\int^{t}_{1}\frac{P(s)}{s^2}ds$:
\begin{equation}
\begin{aligned}
&{\cal E}(\rho,u)(t)=\int_{\R^N}\big(\frac{1}{2}\rho(t,x)|u(t,x)|^{2}+(\Pi(\rho)(t,x)-\Pi(\bar{\rho}))+\kappa|\n\sqrt{\rho}(t,x)|^2\big)dx\\
&+\int^{t}_{0}\int_{\R^N}2\mu\,\rho(t,x)|D u|^{2}(t,x) dt dx \leq \int_{\R^N}\big(\rho_{0}(x)|u_{0}(x)|^{2}+\Pi(\rho_{0})(x)\big)dx.
\end{aligned}
\label{cenergie}
\end{equation}
with $\Pi(\rho)$ defined as follows:
$$\Pi(\rho)=a\big(\rho\ln(\frac{\rho}{\bar{\rho}})+\bar{\rho}-\rho\big)=q(\rho)-q(\bar{\rho})-q'(\bar{\rho})(\rho-\bar{\rho}).$$
Let us observe that in comparison with the compressible Navier-Stokes system the capillary tensor provides additional regularity on the density since the gradient of the square roots of the density is conserved. In fact the capillary tensor makes the system parabolic-parabolic on the density and the velocity.
\begin{remarka}
When $\kappa(\rho)=\frac{\kappa}{\rho}$ we obtain (see the appendix)t:
$${\rm div}K=\kappa{\rm div}(\rho\n \n\ln\rho).$$
\end{remarka}
Let us introduce the unknown $v=u+\mu\n\ln\rho$ (which has been introduced by J\"ungel in \cite{J} in order to prove the existence of global weak solution) such that we have:
$$
\begin{cases}
\begin{aligned}
&\p_t\rho+{\rm div}(\rho v)-\mu\D \rho=0,\\
&\rho\p_t v+\rho u\cdot\n v-\mu{\rm div}(\rho\n u)-\kappa {\rm div}(\rho\n\n\ln\rho)+\n P(\rho)=0.
\end{aligned}
\end{cases}
$$
In particular if $\kappa\geq \mu^2$ we have:
\begin{equation}
\begin{cases}
\begin{aligned}
&\p_t\rho+{\rm div}(\rho v)-\mu\D \rho=0,\\
&\rho\p_t v+\rho u\cdot\n v-\mu{\rm div}(\rho\n v)-(\kappa-\mu^2) {\rm div}(\rho\n\n\ln\rho)+\n P(\rho)=0.
\end{aligned}
\end{cases}
\label{systK}
\end{equation}
Multiplying the previous momentum equation by $v$ when $P(\rho)=a\rho$ we have:
\begin{equation}
\begin{aligned}
&{\cal E}_1(\rho, v)(t)=\int_{\R^N}\big(\frac{1}{2}\rho(t,x)|u(t,x)|^{2}+(\Pi(\rho)(t,x)-\Pi(\bar{\rho}))\big)dx+\int^{t}_{0}\int_{\R^N}\mu\,\rho(t,x)|\n v|^{2}(t,x) dt dx\\
&\hspace{4cm} \leq \int_{\R^N}\big(\rho_{0}(x)|u_{0}(x)|^{2}+(\Pi(\rho_{0})(x)-\Pi(\bar{\rho})\big)dx.
\end{aligned}
\label{bcenergie}
\end{equation}
Combining (\ref{bcenergie}) and (\ref{cenergie}) allows to prove that for any $T>0$ (we refer to \cite{J1}, the result is obtained in the case of the torus, the case of the euclidean space $\R^N$ is a straightforward generalization):
\begin{equation}
\D\sqrt{\rho}\in L^2_T(L^2(\R^N)).
\label{Jungel}
\end{equation}
\begin{remarka}
In the following we shall consider the simple case $\kappa=\mu^2$ and $P(\rho)=a\rho^\gamma$ with $a>0,\gamma\geq 1$. We can observe that when $\kappa=\mu^2$ the previous system can be rewritten as follows:
\begin{equation}
\begin{cases}
\begin{aligned}
&\p_t\rho+{\rm div}(\rho v)-\mu\D \rho=0,\\
&\rho\p_t v+\rho u\cdot\n v-\mu{\rm div}(\rho\n v)+\n P(\rho)=0.
\end{aligned}
\end{cases}
\label{basystK}
\end{equation}
It appears clearly that the density $\rho$ and the effective velocity $v$ verifies parabolic equation. \end{remarka}
\subsection{Motivation and main issue for the Korteweg system}
Now before recalling the main results on the  existence of global weak solutions for compressible Navier Stokes equations and  Korteweg system, we would  like to point out also  an another aspect of the Korteweg system (\ref{13systeme}). Indeed  this system  is also used in a purely theoretical interest consisting in the selection of the physically relevant solutions of the Euler model by a vanishing capillarity-viscosity limit  (especially when the system is not strictly hyperbolic, which is typically the case when the pressure is Van der Waals).  Indeed in this last case at least when $N=1$ it is not possible to apply the classical theory of Lax for the Riemann problem (see  \cite{Lax}) and of Glimm (see \cite{Glimm}) with small $BV$ initial data in order to obtain the existence of global weak-entropy solution (we refer also to the work of Bianchini and Bressan see \cite{BB} for the uniqueness). It explains in particular why it seems important to prove the existence of global strong solution with large initial data for the Korteweg system.\\
In this spirit, we prove recently in \cite{CH2} with  Charve  inspired by DiPerna \cite{Di1} the existence of global strong solution with large initial data of the Korteweg system in one dimension when $\kappa=\mu^2$. To do this we adapt the notion of Riemann invariant to the Korteweg system which allows us to control the $L^\infty$ norm on $v$, next using the maximum principle we can estimate the $L^\infty$ norm on the vacuum $\frac{1}{\rho}$ what is sufficient in order to propagate globally the regularity. In addition we show that this global strong solution converges in the setting of a $\gamma$ law for the pressure ($P(\rho)=a\rho^{\gamma}$, $\gamma>1$) to weak entropy solution of the compressible Euler equations. In particular it justifies that the Korteweg system is suitable for selecting the physical solutions in the case where the Euler system is strictly hyperbolic and when $\kappa=\mu^2$. The problem remains however open for a Van der Waals pressure.\\
Let us give some few words on the choice $\kappa=\mu^{2}$. It gives a specific structure of \textit{effective velocity} (see the system (\ref{basystK}), but is also reasonable on a physic point of view (see \cite{Rohdehdr} for more details). This algebraic relation between $\kappa$ and $\mu^{2}$ corresponds to an intermediary regime, indeed an important research line (see \cite{Rohdehdr}]) is to model the capillarity tensor and to understand how the solutions converges to the Euler system when the capillarity and the viscosity coefficients tends to zero. We want point out here that it exists three different regimes, more precisely if we assume the viscosity coefficient equal to $\e $ with $\e\rightarrow 0$. Then we have the following regimes:
\begin{enumerate}
\item $\kappa<< \e^{2}$, the viscosity dominates so the parabolic effects is primordial. 
\item $\kappa\sim \e^{2}$, intermediary regime.
\item $\kappa>> \e^{2}$, the capillarity dominates so the dispersive effects are predominant.
\end{enumerate}
In particular in the parabolic regime we attend to converge to weak entropy solution, in \cite{CH2} this result is proved in the case of the intermediary regime when $\kappa=\mu^2$. In the last case we probably converge to "dispersive solutions" which generates "dispersive shocks". This problem remains open, it has been studied in the case of Korteweg de Vries equation by Lax and Levermore (see \cite{LL1,LL2,LL3}) using the inverse scattering theory (in this case the solutions converges to a "dispersive solution" of the Burger equation).
\subsection{Existence of global weak solutions for Korteweg system}
We can observe that via the energy inequality (\ref{cenergie}), the density $\sqrt{\rho}$ belongs in
$L^{\infty}(0,\infty,{\dot {H}}^{1}(\R^{N}))$. Hence, in contrast to
the compressible Navier-Stokes system one can easily deal with the
pressure term by involving Sobolev embedding (in other words we have enough compactness on the density). However it appears a new obstruction compared with the compressible Navier-Stokes system which consists in dealing with the quadratic terms in gradient of the density coming from the capillary tensor (see (\ref{divK})).  Bresch, Desjardins and  Lin in \cite{5BDL} got
some stability result for the global weak solutions of the Korteweg model with some
specific viscosity coefficients and capillarity coefficient
$\mu(\rho)=\rho$, $\lambda(\rho)=0$ and $\kappa(\rho)=\kappa$ a constant. To do this they exhibit new entropy inequalities which ensures regularizing effects on the density, roughly speaking they obtain a control of $\D\rho$ in $L^2_T(L^2(\R^N))$ for any $T>0$. It provides enough compactness in order to deal with the capillary tensor. However the global weak solutions of \cite{5BDL} require some
specific test
functions which depend on the solution itself (in other words they obtain the stability of global weak solution for the Korteweg system where the momentum equation is multiplied by $\rho$). This is due to the difficulty to deal with the term $\rho u\otimes u$ (indeed compared withe the case of non degenerate viscosity coefficient we have not directly a control on the gradient of $u$ but only on $\rho\n u$ which is in $L^2_T(L^2(\R^N))$). In \cite{J}, J\"ungel in a very interesting paper obtains by using an effective velocity $v$ the existence of global weak solution when $\kappa(\rho)=\frac{1}{\rho}$ (which corresponds to the quantum compressible Navier-Stokes system) modulo that as in \cite{5BDL} the test functions depends on the density $\rho$. In particular he is the first to introduce this new effective velocity $v$ which allows to simplify the system as we have seen in (\ref{basystK}). It allows him to establish new entropy estimates in the spirit of \cite{5BDL}.\\
In \cite{fH2}, we improve this result  by showing the existence of global weak solution with small initial data in the energy space
for specific choices on the capillary coefficients and with general viscosity coefficient.
Comparing with the results of \cite{5BDL}, we get global weak solutions with general test function $\va\in C^{\infty}_{0}(\R^{N})$ not depending on the density $\rho$. In fact we have extracted of the structure of capillarity term a new energy inequality using fractionary derivative which allows a
gain of derivative on the density $\rho$. Finally in \cite{Hglobal} we extend the result of \cite{J} when $\kappa=\mu^2$ by proving the stability of the weak solution with classical test functions (they do not depend on the solution itself). To do this we involve additional entropy on the effective velocity which provide a gain of integrability on the velocity $v$ (in the spirit of \cite{fMV1}). The existence of global weak solution has been proved recently in an interesting work by P. Antonelli and S. Spirito in \cite{Aa} when $\kappa<\mu^2<\frac{9}{8}\kappa^2$ which corresponds to an intermediary regime. Let us mention in particular that the authors construect global approximate solutions verifying different energy inequalities (the energy entropy, the BD entropy and the Mallet Vasseur Gain of integrability on the effective velocity $v$).
\subsection{Global strong solutions with small initial data for the Korteweg system}
Let us now recall the notion of scaling for the Korteweg's system (\ref{13systeme}). Such an
approach is now classical for incompressible Navier-Stokes equation
and yields local well-posedness (or global well-posedness for small
data) in spaces with minimal regularity.
In our situation we can
easily check that, if $(\rho,u)$ solves (\ref{13systeme}), then
$(\rho_{\lambda},u_{\lambda})$ solves also this system:
$$\rho_{\lambda}(t,x)=\rho(\lambda^{2}t,\lambda x)\,,\,u_{\lambda}(t,x)=\lambda u(\lambda^{2}t,\lambda x)$$
provided the pressure laws $P$ have been changed into
$\lambda^{2}P$.
\begin{definition}
We say that a functional space is critical with respect to the
scaling of the equation if the associated norm is invariant under
the transformation:
$$(\rho,u)\longrightarrow(\rho_{\lambda},u_{\lambda})$$
(up to a constant independent of $\lambda$).
\label{scal}
\end{definition}
This suggests us to choose initial data
$(\rho_{0},u_{0})$ in spaces whose norm is invariant for
all $\lambda>0$ by the transformation
$(\rho_{0},u_{0})\longrightarrow(\rho_{0}(\lambda\cdot),\lambda
u_{0}(\lambda\cdot)).$ A natural candidate is the Besov space (see the section \ref{section2} for some definitions of Besov spaces)
$B^{N/2}_{2,\infty}\times (B^{N/2-1}_{2,\infty})^{N}$,
however since $B^{N/2}_{2,\infty}$ is not included in $L^{\infty}$, we cannot
expect to get a priori $L^{\infty}$ estimate on the density when
$\rho_{0}\in B^{N/2}_{2,\infty}$ (in particular it makes the study of the non linear term delicate since it appears impossible to use composition theorems). Danchin and Desjardins in \cite{fDD} have been the first to obtain the existence of global strong solution with small initial data in the framework of critical Besov spaces choosing $(\rho_0-1,u_0)$ in $B^{\N}_{2,1}\times B^{\N-1}_{2,1}$. Let us mention that this choice allows to control the $L^\infty$ norm on the density since $B^{\N}_{2,1}$ is embedded in $L^\infty$. This last result has been recently improved in \cite{Hprepa} inasmuch as we
can deal with large space of initial data, here  $(\rho_0-1,u_0)$ belong in $(B^{\N}_{2,2}\cap L^\infty)\times B^{\N-1}_{2,2}$. The main difficulty consists in estimating the $L^\infty$ norm of the density, to do this we use a characterization of the Besov space in terms of the semi-group associated to the linearized Korteweg system combined with a maximum principle. Let us also point out that when $\kappa=\mu^2$ we prove the existence of global strong solution with large initial data for the scaling of the equation when $N\geq 2$ on the rotational part. In particular it shows that for some suitable initial data we have existence of global strong solution in dimension $N=2$ however that these initial data are large in the energy space (let us mention that the problem of global strong solution in dimension $N=2$ with large initial data remains open in full generality). The key ingredient of the proof is the notion of quasi-solution, which consists to approximate our solution by an exact solution of the pressure less system of (\ref{13systeme}). We take into account the fact that we have not an exact invariance by scaling because the pressure term and it allows us to consider this term as a small perturbation in high frequencies. We refer also to \cite{global} when we consider the system with friction, we have in this case exact global strong solution with large initial data provided that the velocity is irrotational.  Let us finally cite the work of  Kotschote in \cite{3MK} who showed the existence of strong solution for the isothermal model in bounded domain by using Dore\^a-Venni Theory and $\mathcal{H}^{\infty}$ calculus. We refer also to \cite{3MK1}. In \cite{fH1},  we generalize the results of \cite{fDD} in the case of non isothermal Korteweg system with physical coefficients depending on the density and the temperature.\\
Actually the existence of global strong solution with large initial data in the general case remains a open problem in dimension $N\geq 2$. We would like in the present paper answer to this question when $\kappa(\rho)=\frac{\mu^2}{\rho}$, $\mu(\rho)=2\mu\rho$, $\lambda(\rho)=0$ and $P(\rho)=a\rho^\gamma$ with some conditions on $\gamma\geq 1$. To do this we start by giving a result of strong solution in finite time in critical Besov space, the main issue of our analysis in this first result will be to give an accurate estimate on the time of existence $T^*$. More precisely we shall bounded by below $T^*$ in terms of the norms of the initial data. The next step and the main ingredients of the proof of the existence of global strong solution correspond to provide new energy estimate. In particular we are going to show that we can control $\rho^{\frac{1}{p}}v$ in any $L^\infty_T(L^p(\R^N))$ space with $2\leq p<+\infty$ (in other word we get a gain of integrability on the effective velocity). In the sequel, we shall transfer these new regularity on $\rho^{\frac{1}{p}}v$ on the density $\rho$ via the first equation of (\ref{basystK}) which is parabolic. To realize this program it seems necessary to control the vacuum or the $L^\infty$ norm of $\frac{1}{\rho}$. In order to get such estimates, we will apply method introduced by De Giorgi (see \cite{DG}) to obtain regularity results for elliptic equations with discontinuous diffusion coefficients and 
extended by Ladyzenskaya et al (see \cite{La}) to the parabolic case. One of the main difficulty in order to adapted these technics is that we have no control on $v$ but only on $\rho^{\frac{1}{p}}v$, in order to overcome this difficulty we shall apply a suitable bootstrap argument (indeed $\frac{1}{\rho}$ is bounded in $L^\infty$ norm at least on $(0,T^*)$) which enables us to prove that $\frac{1}{\rho}$ is bounded in reality in $L^\infty_{loc}(L^\infty(\R^N))$. To finish via the accurate estimate on the time existence of a strong solution in terms of the norm of the initial data, we are going to extend our strong solution beyond $T^*$ which will imply that necessary we must have $T^*=+\infty$. Let us state now our mains results.
\section{Main results}
Let us rewrite the system (\ref{systK}) in terms of  $q=\ln(\frac{\rho}{\bar{\rho}})$ and we assume that \begin{equation}
P(\rho)=a\rho^{\gamma},\;\kappa(\rho)=\frac{\mu^2}{\rho}, \;\mu(\rho)=2\mu\rho,\;\lambda(\rho)=0,
\label{hyppothese}
\end{equation}
we have (assuming that the density does not admit vacuum):
\begin{equation}
\begin{cases}
\begin{aligned}
&\p_t q-\mu\D q+u\cdot\n q+{\rm div}v=0,\\
&\p_t v+ u\cdot\n v-\mu\D v-\mu\n q\cdot\n v+a\n P(\rho)=0.
\end{aligned}
\end{cases}
\label{systK1}
\end{equation}
Our first result concerns the existence of strong solution in finite time with in addition an estimate on the time of existence in terms of the initial data.
\begin{theorem}
Let $P$ a general regular pressure. Let $N\geq 2$ and assume that the physical coefficients verify (\ref{hyppothese}). Let  $(q_0,v_0)\in B^{\NN}_{p,1}\times B^{\NN-1}_{p,1}$ with $1\leq p<2N$ and $c$ such that $0<c\leq\rho_0$, then it exists a time $T$ such that system (\ref{systK1}) has a unique solution on $[0,T]$ with:
\begin{equation}
q\in \widetilde{C}_{T}(B^{\NN}_{p,1})\cap L^1_T(B^{\NN+2}_{p,1}),\;\frac{1}{\rho},\rho\in L^{\infty}_{T}(L^\infty(\R^N))\;\;\;\mbox{and}\;\;v\in \widetilde{C}_{T}(B^{\NN-1}_{p,1})\cap L^1_T(B^{\NN+1}_{p,1}).
\label{hypimportant}
\end{equation}
If in addition $(q_0,v_0)$ belongs in $B^{s}_{p,1}\times B^{s-1}_{p,1}$ for any $s>\NN$
then we have:
\begin{equation}
q\in \widetilde{C}_{T}(B^{s}_{p,1})\cap L^1_T(B^{s+2}_{p,1}),\;\frac{1}{\rho},\rho\in L^{\infty}_{T}(L^\infty(\R^N))\;\;\mbox{and}\;\;v\in \widetilde{C}_{T}(B^{s-1}_{p,1})\cap L^1_T
(B^{s+1}_{p,1}).
\label{estimcru33}
\end{equation}
Now we assume that $(q_0,v_0)\in (B^{\NN}_{p,1}\cap B^{\NN+\e}_{p,1})\times (B^{\NN-1}_{p,1}\cap B^{\NN-1+\e}_{p,1})$ with $\e>0$ then it exists $C,C_1,c>0$ such that: 
\begin{equation}
\begin{aligned}
&T\geq  \Inf \big(\frac{ 2(c\mu)^{\frac{2}{\e'}-1}   \e^{\frac{2}{\e'}}}{(8C)^{\frac{2}{\e'}}\|q_0\|_{B^{\NN+\e'}_{p,1}}^{\frac{2}{\e'}}},\frac{ 2(c\mu)^{\frac{2}{\e'}-1}   \e^{\frac{2}{\e'}}}{(8C)^{\frac{2}{\e'}}\|v_0\|_{B^{\NN-1+\e'}_{p,1}}^{\frac{2}{\e'}}}, \frac{C_1}{4},\\
&\hspace{3cm}\frac{1}{16 C^2_1  (\|q_0\|_{B^{\NN}_{p,1}}+\|v_0\|_{B^{\NN-1}_{p,1}})(1+\sqrt{ \|q_0\|_{B^{\NN}_{p,1}}+\|v_0\|_{B^{\NN-1}_{p,1}}})^2)}).
\end{aligned}
\label{superimportant}
\end{equation}
\label{theo1}
\end{theorem}
\begin{remarka}
Let us mention that this theorem is not really new, indeed similar results have been proved in \cite{fDD}. In \cite{Hglobal} it is even possible to obtain strong solution with large class of initial data, it means with $(q_0,v_0)$ in $(B^{\N}_{2,2}\cap L^\infty)\times B^{\N-1}_{2,2}$. Here the main interest of this result provides of the precise estimate (\ref{superimportant}) on the time of existence which will be crucial in order to prove the existence of global strong solution. Let us observe that we need to choose initial data which are slightly surcritical for the scaling of the equation. 
\end{remarka}
Let us define now the maximal time $T^*$ of existence of global strong solution:
$$T^{*}=\sup\{T\in\R\,;\mbox{it exists a solution}\;(q,u)\;\mbox{of the system \ref{systK1} on $[0,T]$ verifying (\ref{hypimportant})} \}$$
\begin{theorem}
Let $N\geq 2$  and assume that the physical coefficients verify (\ref{hyppothese}) with $\gamma\geq1$ when $N=2$, $1\leq\gamma<\frac{7}{3}$ when $N=3$, $\gamma=1$ when $N\geq 4$ . Let $(q_0,v_0)\in (B^{\NN}_{p,1}\cap B^{\NN+\e}_{p,1}) \times (B^{\NN-1}_{p,1}\cap B^{\NN-1+\e}_{p,1})$ with $\e>0$ and with $\frac{N}{1-\e}<p<2N$ and $c$ such that $0<c\leq\rho_0$. In addition we assume
$v_0\in L^\infty$ and $(\rho,v_0)$ are of finite energy which implies:
$$
\begin{aligned}
&{\cal E}(\rho_0,u_0)<+\infty,\\
&{\cal E}_1(\rho_0,v_0)<+\infty.
\end{aligned}
$$
then  we have $T^{*}=+\infty$. This global solution is unique.
\label{theo2}
\end{theorem}
\begin{remarka}
The above assumption on $v_0$ ensures that $\rho_0^{\frac{1}{p_1}}v_0$ is uniformly bounded  in $p_1$ in $L^{p_1}$ for any $2\leq p_1<+\infty$. It will be important in the sequel in order to get a gain of integrability on the velocity $v$.
\end{remarka}
\begin{remarka}
This is the first result up our knowledge of global strong solution with large initial data in full generality when $N\geq 2$ for the Korteweg system with the specific choice on the physical coefficients (\ref{hyppothese}). Let us mention that this particular choice on the coefficients enables us to obtain gain of integrability on the effective velocity $v$ which is crucial in the proof of our result.
\end{remarka}
\begin{remarka}
This result allows also to prove the existence of global weak solution with initial data in the energy space with this particular choice on the physical coefficients. Indeed in \cite{Hglobal} we prove the stability of the global weak solution, it was however difficult to construct approximate global weak solution which verify uniformly the different entropy inequalities. This is done now using the previous theorem.
\end{remarka}
The paper is structured in the following way: in section \ref{section2} we recall some important results on the Littlewood-Paley theory and the notion of Besov spaces. In section \ref{section3}, we prove the theorem \ref{theo1}. In the section \ref{section4} we show uniform estimate in time on the solution $(\rho, v)$. In the last section \ref{section5} we end up with the proof of the theorem \ref{theo2} when $P(\rho)=a\rho$ with $a>0$. In section \ref{section6} we prove the theorem \ref{theo2}  for general pressure. To conclude we postpone in the appendix the derivation of the model (\ref{systK}).

\section{Littlewood-Paley theory and Besov spaces}
\label{section2}
Throughout the paper, $C$ stands for a constant whose exact meaning depends on the context. The notation $A\lesssim B$ means
that $A\leq CB$.\\
As usual, the Fourier transform of $u$ with respect to the space variable will be denoted by $\mathcal{F}(u)$ or $\hat{u}$. 
In this section we will state classical definitions and properties concerning the homogeneous dyadic decomposition with respect to the Fourier variable. We will recall some classical results and we refer to \cite{BCD} (Chapter 2) for proofs (and more general properties).

To build the Littlewood-Paley decomposition, we need to fix a smooth radial function $\chi$ supported in (for example) the ball $B(0,\frac{4}{3})$, equal to 1 in a neighborhood of $B(0,\frac{3}{4})$ and such that $r\mapsto \chi(r.e_r)$ is nonincreasing over $\R_+$. So that if we define $\varphi(\xi)=\chi(\xi/2)-\chi(\xi)$, then $\varphi$ is compactly supported in the annulus $\{\xi\in \R^d, \frac{3}{4}\leq |\xi|\leq \frac{8}{3}\}$ and we have that,
\begin{equation}
 \forall \xi\in \R^d\setminus\{0\}, \quad \sum_{l\in\Z} \varphi(2^{-l}\xi)=1.
\label{LPxi}
\end{equation}
Then we can define the \textit{dyadic blocks} $(\ddl)_{l\in \Z}$ by $\ddl:= \varphi(2^{-l}D)$ (that is $\hat{\ddl u}=\varphi(2^{-l}\xi)\hat{u}(\xi)$) so that, formally, we have
\begin{equation}
u=\Sum_l \ddl u
\label{LPsomme} 
\end{equation}
As (\ref{LPxi}) is satisfied for $\xi\neq 0$, the previous formal equality holds true for tempered distributions \textit{modulo polynomials}. A way to avoid working modulo polynomials is to consider the set $\cS_h'$ of tempered distributions $u$ such that
$$
\lim_{l\rightarrow -\infty} \|\dot{S}_l u\|_{L^\infty}=0,
$$
where $\dot{S}_l$ stands for the low frequency cut-off defined by $\dot{S}_l:= \chi(2^{-l}D)$. If $u\in \cS_h'$, (\ref{LPsomme}) is true and we can write that $\dot{S}_l u=\Sum_{q\leq l-1} \ddq u$. We can now define the homogeneous Besov spaces used in this article:
\begin{definition}
\label{LPbesov}
 For $s\in\R$ and  
$1\leq p,r\leq\infty,$ we set
$$
\|u\|_{B^s_{p,r}}:=\bigg(\sum_{l} 2^{rls}
\|\Delta_l  u\|^r_{L^p}\bigg)^{\frac{1}{r}}\ \text{ if }\ r<\infty
\quad\text{and}\quad
\|u\|_{ B^s_{p,\infty}}:=\sup_{l} 2^{ls}
\|\Delta_l  u\|_{L^p}.
$$
We then define the space $ B^s_{p,r}$ as the subset of  distributions $u\in {\cS}'_h$ such that $\|u\|_{ B^s_{p,r}}$ is finite.
\end{definition}
Once more, we refer to \cite{BCD} (chapter $2$) for properties of the inhomogeneous and homogeneous Besov spaces. Among these properties, let us mention:
\begin{itemize}
\item for any $p\in[1,\infty]$ we have the following chain of continuous embeddings:
$$
 B^0_{p,1}\hookrightarrow L^p\hookrightarrow  B^0_{p,\infty};
$$
\item if $p<\infty$ then 
  $B^{\frac dp}_{p,1}$ is an algebra continuously embedded in the set of continuous 
  functions decaying to $0$ at infinity;
    \item for any  smooth homogeneous  of degree $m$ function $F$ on $\R^d\setminus\{0\}$
the operator $F(D)$ defined by $F(D)u=\mathcal{F}^{-1}\Big(F(\cdot)\mathcal{F}(u)(\cdot)\Big)$ maps  $ B^s_{p,r}$ in $ B^{s-m}_{p,r}.$ This implies that the gradient operator maps $ B^s_{p,r}$ in $B^{s-1}_{p,r}.$  
  \end{itemize}
We refer to \cite{BCD} (lemma 2.1) for the Bernstein lemma (describing how derivatives act on spectrally localized functions), that entails the following embedding result:
\begin{proposition}\label{LP:embed}
\sl{For all $s\in\R,$ $1\leq p_1\leq p_2\leq\infty$ and $1\leq r_1\leq r_2\leq\infty,$
  the space $ B^{s}_{p_1,r_1}$ is continuously embedded in 
  the space $ B^{s-d(\frac1{p_1}-\frac1{p_2})}_{p_2,r_2}.$}
\end{proposition}
Then we have:
$$
 B^\fdp_{p,1}\hookrightarrow  B^0_{\infty,1}\hookrightarrow L^\infty.
$$
In this paper, we shall mainly work with functions or distributions depending on both the time variable $t$ and the space variable $x.$ We shall denote by $\cC(I;X)$ the set of continuous functions on $I$ with values in $X.$ For $p\in[1,\infty]$, the notation $L^p(I;X)$ stands for the set of measurable functions on  $I$ with values in $X$ such that $t\mapsto \|f(t)\|_X$ belongs to $L^p(I)$.

In the case where $I=[0,T],$  the space $L^p([0,T];X)$ (resp. $\cC([0,T];X)$) will also be denoted by $L_T^p X$ (resp. $\cC_T X$). Finally, if $I=\R^+$ we shall alternately use the notation $L^p X.$

The Littlewood-Paley decomposition enables us to work with spectrally localized (hence smooth) functions rather than with rough objects. We naturally obtain bounds for each dyadic block in spaces of type $L^\rho_T L^p.$  Going from those type of bounds to estimates in  $L^\rho_T \dot B^s_{p,r}$ requires to perform a summation in $\ell^r(\Z).$ When doing so however, we \emph{do not} bound the $L^\rho_T \dot B^s_{p,r}$ norm for the time integration has been performed \emph{before} the $\ell^r$ summation.
This leads to the following notation:

\begin{definition}\label{d:espacestilde}
For $T>0,$ $s\in\R$ and  $1\leq r,\sigma\leq\infty,$
 we set
$$
\|u\|_{\tilde L_T^\sigma  B^s_{p,r}}:=
\bigl\Vert2^{js}\|\ddq u\|_{L_T^\sigma L^p}\bigr\Vert_{\ell^r(\Z)}.
$$
\end{definition}
One can then define the space $\tilde L^\sigma_T \dot B^s_{p,r}$ as the set of  tempered distributions $u$ over $(0,T)\times \R^d$ such that $\lim_{q\rightarrow-\infty}\dot S_q u=0$ in $L^\sigma([0,T];L^\infty(\R^d))$ and $\|u\|_{\tilde L_T^\sigma  B^s_{p,r}}<\infty.$ The letter $T$ is omitted for functions defined over $\R^+.$ 
The spaces $\tilde L^\sigma_T  B^s_{p,r}$ may be compared with the spaces  $L_T^\sigma \dot B^s_{p,r}$ through the Minkowski inequality: we have
$$
\|u\|_{\tilde L_T^\sigma  B^s_{p,r}}
\leq\|u\|_{L_T^\sigma  B^s_{p,r}}\ \text{ if }\ r\geq\sigma\quad\hbox{and}\quad
\|u\|_{\tilde L_T^\sigma  B^s_{p,r}}\geq
\|u\|_{L_T^\sigma  B^s_{p,r}}\ \text{ if }\ r\leq\sigma.
$$
All the properties of continuity for the product and composition which are true in Besov spaces remain true in the above  spaces. The time exponent just behaves according to H\"older's inequality. 
\medbreak
Let us now recall a few nonlinear estimates in Besov spaces. Formally, any product of two distributions $u$ and $v$ may be decomposed into 
\begin{equation}\label{eq:bony}
uv=T_uv+T_vu+R(u,v), \mbox{ where}
\end{equation}
$$
T_uv:=\sum_l\dot S_{l-1}u\ddl v,\quad
T_vu:=\sum_l \dot S_{l-1}v\ddl u\ \hbox{ and }\ 
R(u,v):=\sum_l\sum_{|l'-l|\leq1}\ddl u\,\dot\Delta_{l'}v.
$$
The above operator $T$ is called ``paraproduct'' whereas $R$ is called ``remainder''. The decomposition \eqref{eq:bony} has been introduced by Bony in \cite{BJM}.

In this article we will frequently use the following estimates (we refer to \cite{BCD} section 2.6):.
\begin{proposition}
Under the same assumptions there exists a constant $C>0$ such that if $1/p_1+1/p_2=1/p$, and $1/r_1+1/r_2=1/r$:
$$
\|\dot{T}_u v\|_{B_{2,1}^s}\leq C \|u\|_{L^\infty} \|v\|_{B_{2,1}^s},
$$
$$\|\dot{T}_u v\|_{B_{p,r}^{s+t}}\leq C\|u\|_{B_{p_1,r_1}^t} \|v\|_{\dot{B}_{p_2,r_2}^s} \quad (t<0),
$$
\begin{equation}
 \|\dot{R}(u,v)\|_{B_{p,r}^{s_1+s_2-\fd}} \leq C\|u\|_{B_{p_1,r_1}^{s_1}} \|v\|_{B_{p_2,r_2}^{s_2}} \quad (s_1+s_2>0).
\end{equation}
\label{produit}
\end{proposition}
Let us now turn to the composition estimates. We refer for example to \cite{BCD} (Theorem $2.59$, corollary $2.63$)):
\begin{proposition}
\sl{\begin{enumerate}
 \item Let $s>0$, $u\in B_{p,1}^s\cap L^{\infty}$ and $F\in W_{loc}^{[s]+2, \infty}(\R^d)$ such that $F(0)=0$. Then $F(u)\in B_{p,1}^s$ and there exists a function of one variable $C_0$ only depending on $s$, $p$, $d$ and $F$ such that
$$
\|F(u)\|_{B_{p,1}^s}\leq C_0(\|u\|_{L^\infty})\|u\|_{B_{p,1}^s}.
$$
\item If $u$ and $v\in B_{p,1}^\fd$ and if $v-u\in B_{p,1}^s$ for $s\in]-\fd, \fd]$ and $G\in W_{loc}^{[s]+3, \infty}(\R^d)$, then $G(v)-G(u)$ belongs to $B_{p,1}^s$ and there exists a function of two variables $C$ only depending on $s$, $d$ and $G$ such that
$$
\|G(v)-G(u)\|_{B_{p,1}^s}\leq C(\|u\|_{L^\infty}, \|v\|_{L^\infty})\left(|G'(0)| +\|u\|_{B_{p,1}^\fd} +\|v\|_{B_{p,1}^\fd}\right) \|v-u\|_{B_{p,1}^s}.
$$
\end{enumerate}}
\label{composition}
\end{proposition}
Let us now recall a result of interpolation which explains the link between the space $B^{s}_{p,1}$ and the space $B^{s}_{p,\infty}$ (see
\cite{BCD} sections $2.11$ and $10.2.4$):
\begin{proposition}\sl{
\label{interpolation}
There exists a constant $C$ such that for all $s\in\R$, $\e>0$, $\sigma\geq 1$ and
$1\leq p<+\infty$,
$$\|u\|_{\widetilde{L}_{T}^{\sigma}(B^{s}_{p,1})}\leq C\frac{1+\e}{\e}\|u\|_{\widetilde{L}_{T}^{\sigma}(B^{s}_{p,\infty})}
\log\biggl(e+\frac{\|u\|_{\widetilde{L}_{T}^{\sigma}(B^{s-\e}_{p,\infty})}+ \|u\|_{\widetilde{L}_{T}^{\sigma}(B^{s+\e}_{p,\infty})}}
{\|u\|_{\widetilde{L}_{T}^{\sigma}(B^{s}_{p,\infty})}}\biggl).$$ \label{5Yudov}}
\end{proposition}
\subsection*{Parabolic equations}
Let us end this section by recalling the following estimates for the heat equation:
\begin{proposition}\sl{
\label{chaleur} Let $s\in\R$, $(p,r)\in[1,+\infty]^{2}$ and
$1\leq\rho_{2}\leq\rho_{1}\leq+\infty$. Assume that $u_{0}\in B^{s}_{p,r}$ and $f\in\widetilde{L}^{\rho_{2}}_{T}
(B^{s-2+2/\rho_{2}}_{p,r})$.
Let u be a solution of:
$$
\begin{cases}
\begin{aligned}
&\p_{t}u-\mu\D u=f\\
&u_{/t=0}=u_{0},\\
\end{aligned}
\end{cases}
$$  
where $\mu>0$. Then there exists $C>0$ depending only on $N,\mu,\rho_{1}$ and
$\rho_{2}$ such that:
$$\|u\|_{\widetilde{L}^{\rho_{1}}_{T}(B^{s+2/\rho_{1}}_{p,r})}\leq C\big(
 \|u_{0}\|_{B^{s}_{p,r}}+\|f\|_{\widetilde{L}^{\rho_{2}}_{T}
 (B^{s-2+2/\rho_{2}}_{p,r})}\big)\,.$$
 If in addition $r$ is finite then $u$ belongs to $C([0,T],B^{s}_{p,r})$. We have in fact for $\rho'_1=(1+\frac{1}{\rho_1}-\frac{1}{\rho_2})^{-1}$:
  \begin{equation}
  \begin{aligned}
&\|u_L\|_{\widetilde{L}^{\rho_1}_{T}(B^{s+\frac{2}{\rho_1}}_{p,1})}\leq C\big(\sum_{q\in\mathbb{Z}}2^{qs}\|\D_q u_0\|_{L^p}\big(\frac{1-e^{-c\mu T2^{2q}\rho_1}}{c\mu\rho_1})^{\frac{1}{\rho_1}}\\
&\hspace{2cm}+\sum_{q\in\mathbb{Z}}2^{q(s-2+\frac{2}{\rho_2})}\|\D_q f\|_{L^{\rho_2}_T(L^p)}\big(\frac{1-e^{-c\mu T2^{2q}\rho_1}}{c\mu\rho'_1})^{\frac{1}{\rho'_1}}\big)
\end{aligned}
\label{estimtempsimp}
\end{equation}}
\end{proposition}
\section{Proof of the theorem \ref{theo1}}
\label{section3}
\subsection{Existence of local solutions for system (\ref{13systeme})}
\label{sub32}
We now are going to prove the existence of strong solutions in finite time with large initial data verifying the hypothesis of theorem \ref{theo1} for the system (\ref{systK1}). In addition we will provide estimate on the time of existence. We assume now that $(q_{0},v_{0})$ belongs to $B^{\NN}_{p,1}\times B^{\NN-1}_{p,1}$ with $1\leq p<2N$. 
\subsubsection*{Existence of solutions}
For simplicity we just consider the case $P(\rho)=a\rho$, the general cases a straightforward adaptation.  Let us recall the form of the system (\ref{systK1}):
\begin{equation}
\begin{cases}
\begin{aligned}
&\p_t q-\mu\D q+{\rm div}v=-v\cdot\n q+\mu|\n q|^2,\\
&\p_t v+u\cdot\n v-\mu\D v-\mu\n q\cdot\n v+a\n q=0.
\end{aligned}
\end{cases}
\end{equation}
The existence part of the theorem is proved by an iterative method. We define a sequence $(q^{n},u^{n})$ as follows:
$$q^{n}=q_{L}+\bar{q}^{n},\;u^{n}=u_{L}+\bar{v}^{n},$$
where $(q_{L},u_{L})$ stands for the  solution of:
\begin{equation}
\begin{cases}
\p_{t}q_{L}+{\rm div}v_{L}-\mu\D q_L=0,\\
\p_{t}v_{L}-\mu \D v_{L}=0,
\end{cases}
\label{lineaire}
\end{equation}
supplemented with initial data:
$$q_{L}(0)=q_{0}\;,\;v_{L}(0)=v_{0}.$$
Using the proposition \ref{chaleur}, we obtain the following estimates on $(q_{L},u_{L})$ for all $T>0$:
$$q_{L}\in\widetilde{C}([0,T],B^{\NN}_{p,1})\cap\widetilde{L}^{1}_{T}(
B^{\N+2}_{p,1})\;\;\mbox{and}\;\;v_{L}\in\widetilde{C}([0,T],B^{\NN-1}_{p,1})\cap
\widetilde{L}^{1}_{T}
(B^{\N+1}_{p,1}).$$
Setting $(\bar{q}^{0},\bar{v}^{0})=(0,0)$ we now define $(\bar{q}_{n},\bar{v}_{n})$ as the solution of the following system:
$$
\begin{cases}
\begin{aligned}
&\p_{t}\bar{q}^{n}+{\rm div}(\bar{v}^{n})-\mu\D \bar{q}^{n}=F_{n-1},\\
& \p_{t}\bar{v}_{n}-\mu\D\bar{v}_{n}=G_{n-1},\\
&(\bar{q}_{n},\bar{u}_{n})_{t=0}=(0,0),
\end{aligned}
\end{cases}
\leqno{(N_{1})}
$$
where:
$$
\begin{aligned}
F_{n-1}=&-v_{n-1}\cdot\n q^{n-1}+\mu|\n q^{n-1}|^2,\\
G_{n-1}=&-u^{n-1}\cdot\n v^{n-1}+\mu\n q^{n-1}\cdot\n v^{n-1}-a\n q^{n-1}.
\end{aligned}
$$
\subsubsection*{1) First Step , Uniform Bound}
 Let $\e$ be a small
parameter and  choose $T$ small enough  such that according to the proposition \ref{chaleur} we have:
$$
\begin{aligned}
&\|v_{L}\|_{\widetilde{L}^{1}_{T}(B^{\NN+1}_{p,1})}+\|q_{L}\|_{\widetilde{L}^{1}_{T}(B^{\NN+2}_{p,1})}\leq 2\e,\\
&\|v_{L}\|_{\widetilde{L}^{\infty}_{T}(B^{\NN-1}_{p,1})}+\|q_{L}\|_{\widetilde{L}^{\infty}_{T}(B^{\NN}_{p,1})}\leq
C A_{0},
\end{aligned}
\leqno{({\cal{H}}_{\e})}
$$
with $A_{0}=\|q_{0}\|_{B^{\NN}_{p,1}}+\|v_{0}\|_{B^{\NN-1}_{p,1}}$. We are going to show by induction that:
$$\|(\bar{q}^{n},\bar{u}^{n})\|_{F_{T}}\leq\sqrt{\e}.\leqno{({\cal{P}}_{n})}$$
for $\e$ small enough with:
$$F_{T}=\big(\widetilde{C}([0,T],B^{\NN}_{p,1})\cap\widetilde{L}^{1}_{T}(
B^{\NN+2}_{p,1})\big)\times\big(\widetilde{C}([0,T],B^{\NN-1}_{p,1})\cap
\widetilde{L}^{1}_{T}
(B^{\NN+1}_{p,1})\big)^{N}.$$
As $(\bar{q}^{0},\bar{u}^{0})=(0,0)$ the result is true for $n=0$. We now suppose $({\cal P}_{n-1})$ (with $n\geq 1$) true and we are going to show  $({\cal P}_{n})$.
Applying proposition \ref{chaleur}  we have:
\begin{equation}
\begin{aligned}
&\|(\bar{q}^{n},\bar{v}^{n})\|_{F_{T}}\leq C\|(\n
F_{n-1},G_{n-1})\|_{\widetilde{L}^{1}_{T}(B^{\NN-1}_{p,1})}.
\end{aligned}
\label{fi11a}
\end{equation}
Bounding the right-hand side of (\ref{fi11a}) may be
done by applying
proposition  \ref{produit} and proposition  \ref{composition}. We begin with treating the case of $\|F_{n-1}\|_{\widetilde{L}^{1}_{T}(B^{\NN}_{p,1})}$, let us recall that:
\begin{equation}
\begin{aligned}
F_{n-1}=&-v_L \cdot\n q_L-v_L \cdot\n \bar{q}^{n-1}-\bar{v}^{n-1}\cdot\n q_L-\bar{v}^{n-1}\cdot\n \bar{q}^{n-1}+\mu|\n q_L|^2\\
&+2\mu\n q_L\cdot\n \bar{q}^{n-1}+\mu |\n \bar{q}^{n-1}|^2.
\end{aligned}
\label{Fn1}
\end{equation}
We are going to bound each term of (\ref{Fn1}), we have then for $C>0$ large enough:
\begin{equation}
\begin{aligned}
\|v_L \cdot\n q_L   \|_{\widetilde{L}^{1}_{T}(B^{\NN}_{p,1})}&\lesssim \|v_L \|_{\widetilde{L}^{2}_{T}(B^{\N}_{p,1})}\| q_L\|_{\widetilde{L}^{2}_{T}(B^{\NN+1}_{p,1})},\\
&\leq  CA_0\e.
\end{aligned}
\label{F1}
\end{equation}
Similarly we obtain for $C$ large enough and $1\leq p<2N$
\begin{equation}
\begin{aligned}
\|v_{L}\cdot\n \bar{q}^{n-1}\|_{\widetilde{L}^{1}_{T}(B^{\NN}_{p,1})}&\leq C
\|\n \bar{q}^{n-1}\|_{\widetilde{L}^{2}_{T}(B^{\NN}_{p,1})}\|v_{L}\|_{\widetilde{L}^{2}_{T}(B^{\NN}_{p,1})},\\
&\leq C\sqrt{A_0} \e.
\end{aligned}
\label{F2}
\end{equation}
\begin{equation}
\begin{aligned}
\|\bar{v}^{n-1}\cdot\n q_{L}\|_{\widetilde{L}^{1}_{T}(B^{\NN}_{p,1})}&\lesssim \|\bar{v}^{n-1}\|_{\widetilde{L}^{2}_{T}(B^{\NN}_{p,1})}\|\n q_{L}\|_{\widetilde{L}^{2}_{T}(B^{\NN}_{p,1})},\\
&\leq C\e\sqrt{A_0}.
\end{aligned}
\label{F3}
\end{equation}
and:
\begin{equation}
\begin{aligned}
\|\bar{v}^{n-1}\cdot\n \bar{q}^{n-1}\|_{\widetilde{L}^{1}_{T}(B^{\NN}_{p,1})}&\lesssim \|\bar{v}^{n-1}\|_{\widetilde{L}^{2}_{T}(B^{\NN}_{p,1})}\|\n \bar{q}^{n-1}\|_{\widetilde{L}^{2}_{T}(B^{\NN}_{p,1})},\\
&\leq C\e.
\end{aligned}
\label{F4}
\end{equation}
Similarly we have:
\begin{equation}
\begin{aligned}
&\| \,|\n q^{n-1}|^2\|_{\widetilde{L}^{1}_{T}(B^{\NN}_{p,1})}\lesssim \| \n q^{n-1}\|^2_{\widetilde{L}^{2}_{T}(B^{\NN}_{p,1})}\lesssim (\sqrt{\e}(1+\sqrt{CA_0}))^2
\end{aligned}
\label{F4}
\end{equation}
By using the previous inequalities (\ref{F1}), (\ref{F2}), (\ref{F3}), (\ref{F4}), $({\cal P}_{n-1})$ and by interpolation, we obtain that for $C>0$ large enough:
\begin{equation}
\|F_{n}\|_{\widetilde{L}^{1}_{T}(B^{\NN}_{p,1})}\leq C\e\big(CA_0+1+2\sqrt{CA_0}\big).
\label{Fn}
\end{equation}
Next we want to control  $\|G_{n}\|_{\widetilde{L}^{1}(B^{\NN-1}_{p,1})}$. According to
propositions \ref{produit} we have when $2\NN-1>0$ (which is equivalent to $1\leq p<2N$) for $C>0$ large enough:
\begin{equation}
\begin{aligned}
&\|v_L\cdot \n
v_L\|_{\widetilde{L}^{1}_{T}(B^{\NN-1}_{p,1})}\lesssim\|v_L\|_{\widetilde{L}^{2}_{T}(B^{\NN}_{p,1})}\|\n v_L\|_{\widetilde{L}^{2}_{T}(B^{\NN-1}_{p,1})} 
\leq CA_0\e,\\[2mm]
&\|v_L\cdot \n
\bar{v}^{n-1}\|_{\widetilde{L}^{1}_{T}(B^{\NN-1}_{p,1})}\lesssim
\|v_L\|_{\widetilde{L}^{2}_{T}(B^{\NN}_{p,1})}\|\n \bar{v}^{n-1}\|_{\widetilde{L}^{2}_{T}(B^{\NN-1}_{p,1})}\leq 
C \sqrt{ A_0}\e,
\end{aligned}
\label{aGn1}
\end{equation}
\begin{equation}
\begin{aligned}
&\|\bar{v}^{n-1}\cdot \n
v_L\|_{\widetilde{L}^{1}_{T}(B^{\NN-1}_{p,1})}\leq 
C\sqrt{ A_0}\e,\\[2mm]
&\|\bar{v}^{n-1}\cdot \n
\bar{v}^{n-1}\|_{\widetilde{L}^{1}_{T}(B^{\NN-1}_{p,1})}\leq C\e.\\
&\|\n q^{n-1}\cdot \n v^{n-1}\|_{\widetilde{L}^{1}_{T}(B^{\NN-1}_{p,1})}\lesssim \|\n q^{n-1}\|_{\widetilde{L}^{2}_{T}(B^{\NN}_{p,1})}  \| \n v^{n-1}\|_{\widetilde{L}^{2}_{T}(B^{\NN-1}_{p,1})}\leq 
C(\sqrt{\e}+\sqrt{A_0}\sqrt{\e})^2\\
&\hspace{4,5cm}\leq C \e(1+2\sqrt{A_0}+A_0) \\[2mm]
&\|\n q^{n-1}\|_{\widetilde{L}^{1}_{T}(B^{\NN-1}_{p,1})}\leq T\|q^{n-1}\|_{\widetilde{L}^{\infty}_{T}(B^{\NN}_{p,1})}\leq T(CA_0+\sqrt{\e}).
\end{aligned}
\label{aGn}
\end{equation}
Using (\ref{fi11a}), (\ref{Fn}), (\ref{aGn1}), (\ref{aGn})  we obtain for a certain $C_1>0$ large enough and with $\e<1$:
$$
\begin{aligned}
\|(\bar{q}^{n},\bar{u}^{n})\|_{F_{T}}&\leq C_1\e(1+2\sqrt{A_0}+A_0)+C_1T(A_0+\sqrt{\e}).
\end{aligned}
$$
By choosing $T$ and $\e$ small enough  the property $({\cal{P}}_{n})$ is verified, so we
have shown by induction that $(q^{n},u^{n})$ is bounded
in $F_{T}$. To do this we are going to take:
\begin{equation}
\begin{aligned}
&\sqrt{\e}\leq\frac{1}{4C_1(1+2\sqrt{A_0}+A_0)}\;\;\mbox{and}\;\;T\leq\min(\frac{\sqrt{\e}}{4C_1 A_0},\frac{C_1}{4}).
\end{aligned}
\end{equation}
It implies that $T$ must verify the hypothesis $({\cal H}_\e)$ and:
\begin{equation}
\begin{aligned}
&T\leq \min(\frac{1}{16 C^2_1  (\|q_0\|_{B^{\NN}_{p,1}}+\|v_0\|_{B^{\NN-1}_{p,1}})(1+\sqrt{ (\|q_0\|_{B^{\NN}_{p,1}}+\|v_0\|_{B^{\NN-1}_{p,1}})})^2)},\frac{C_1}{4}).
\end{aligned}
\label{impborneT}
\end{equation}
We will show in the next section how to explicit the condition $({\cal H}_\e)$ in the case of slightly surcritical assumption on the initial data.
\subsubsection*{Second Step: Convergence of the
sequence}
 We will show
that $(q^{n},v^{n})$ is a Cauchy sequence in the Banach
space $F_{T}$, hence converges to some
$(q,v)\in F_{T}$. Let:
$$\delta q^{n}=q^{n+1}-q^{n},\;\delta v^{n}=v^{n+1}-v^{n}.$$
The system verified by $(\de q^{n},\de u^{n})$ reads:
$$
\begin{cases}
\begin{aligned}
&\p_{t}\delta q^{n}+{\rm div}\delta v^{n}-\mu\D \delta q^n=F_{n}-F_{n-1},\\
&\p_{t}\delta v^{n}-\mu\D\delta v^{n}=G_{n}-G_{n-1},\\
&\delta q^{n}(0)=0\;,\;\delta v^{n}(0)=0,
\end{aligned}
\end{cases}
$$
Applying propositions \ref{chaleur}, we obtain:
\begin{equation}
\begin{aligned}
\|(\de q^{n},\de u^{n})\|_{F_{T}}\leq\;&
C(\|F_{n}-F_{n-1}\|_{\widetilde{L}^{1}_{T}(B^{\NN}_{p,1})}+\|G_{n}-G_{n-1}\|_{\widetilde{L}^{1}_{T}(B^{\NN-1}_{p,1})}).
\end{aligned}
\label{cauchy}
\end{equation}
Tedious calculus ensure that:
$$
\begin{aligned}
&F_{n}-F_{n-1}=-\delta v^{n-1}\cdot\n q^{n}-v^{n-1}\cdot\n\delta q^{n-1}+\mu\n\delta q^{n-1}\cdot\n q^n+\mu\n\delta q^{n-1}\cdot\n q^{n-1},\\[2mm]
&G_{n}-G_{n-1}=-u^{n}\cdot\n\delta v^{n-1}-\delta u^{n-1}\cdot\n v^{n-1}+\mu\n q^n\cdot\n\delta v^{n-1}+\mu\delta q^{n-1}\cdot\n v^{n-1}-a\n\delta q^{n-1}.
\end{aligned}
$$
It remains only to estimate the terms on the right hand side of (\ref{cauchy}) by using the same type of estimates than in the previous section and the property $({\cal P}_{n})$. More precisely we have via the proposition \ref{produit} and $({\cal P}_{n})$, it exists $C>0$ such that:
\begin{equation}
\begin{aligned}
&\|F_{n}-F_{n-1}\|_{\widetilde{L}^{1}_{T}(B^{\NN}_{p,1})}\lesssim \|\delta u^{n-1}\|_{\widetilde{L}^{\frac{4}{3}}_{T}(B^{\N+\frac{1}{2}}_{p,1})}\|\n q^{n}\|_{\widetilde{L}^{4}_{T}(B^{\NN-\frac{1}{2}}_{p,1})}\\
&\hspace{0,5cm}+ \|\delta u^{n-1}\|_{\widetilde{L}^{\infty}_{T}(B^{\NN-1}_{p,1})}\|\n q^{n}\|_{\widetilde{L}^{1}_{T}(B^{\NN+1}_{p,1})}+\|u^{n-1}\|_{\widetilde{L}^{1}_{T}(B^{\NN+1}_{p,1})}\|\n \delta q^{n-1}\|_{\widetilde{L}^{\infty}_{T}(B^{\NN-1}_{p,1})}\\
&\hspace{6,5cm}+\|\n \delta q^{n-1}\|_{\widetilde{L}^{\frac{4}{3}}_{T}(B^{\NN+\frac{1}{2}}_{p,1})}\|u^{n-1}\|_{\widetilde{L}^{4}_{T}(B^{\NN-\frac{1}{2}}_{p,1})}\\[1,5mm]
&\leq C (A_{0}^{\frac{3}{4}}\e^{\frac{1}{4}}+A_{0}^{\frac{1}{4}}\e^{\frac{3}{4}}+\sqrt{\e}+\e)\|(\delta q^{n-1},\delta v^{n-1})\|_{F_{T}}.
\end{aligned}
\label{cFn}
\end{equation}
In a similar way we show that it exists $C>0$ large enough such that:
\begin{equation}
\begin{aligned}
&\|G_{n}-G_{n-1}\|_{\widetilde{L}^{1}_{T}(B^{\NN-1}_{p,1})}\leq C (A_{0}^{\frac{3}{4}}\e^{\frac{1}{4}}+A_{0}^{\frac{1}{4}}\e^{\frac{3}{4}}+\sqrt{\e}+\e+T)\|(\delta q^{n-1},\delta v^{n-1})\|_{F_{T}}.
\end{aligned}
\label{cGn}
\end{equation}
By combining (\ref{cauchy}), (\ref{cFn}) and (\ref{cGn}) , we get for $C>0$ large enough:
$$\|(\de q^{n},\de v^{n})\|_{F_{T}}\leq C (A_{0}^{\frac{3}{4}}\e^{\frac{1}{4}}+A_{0}^{\frac{1}{4}}\e^{\frac{3}{4}}+\sqrt{\e}+\e+T)\|(\delta q^{n-1},\delta v^{n-1})\|_{F_{T}}.$$ 
It implies that choosing $\e$ and $T$ small enough 
$(q^{n},v^{n})$ is a Cauchy sequence in $F_{T}$ which is a Banach space. It provides that $(q^{n},v^{n})$ converges to
$(q,v)$ is in $F_{T}$. The verification that the limit $(q,v)$ is solution of (\ref{systK1}) in the sense of distributions is a straightforward application of proposition \ref{produit}.
\subsubsection*{Third step: Uniqueness}
Now, we are going to prove the uniqueness of the solution in $F_{T}$. Suppose that $(q_{1},v_{1})$ and $(q_{2},v_{2})$ are solutions with the same initial conditions and belonging in $F_{T}$ where
$(q_{1},v_{1})$ corresponds to the previous
solution. We set:
$$\de q=q_{2}-q_{1}\;\;\;\mbox{and}\;\;\;\de v=v_{2}-v_{1}.$$
We deduce that $(\de q,\de v)$ satisfy the following system:
$$
\begin{cases}
\begin{aligned}
&\p_{t}\delta q+{\rm div}\delta
v-\mu\D\delta q=F_{2}-F_{1},\\
&\p_{t}\delta v-\mu\D\delta v=G_{1}-G_{2},\\
&\delta q(0)=0\;,\;\delta u(0)=0.
\end{aligned}
\end{cases}
$$
We now apply proposition \ref{chaleur} to the previous system, and by using the same type of estimates than in the previous part, we
show that:
$$
\begin{aligned}
&\|(\de q,\de v)\|_{\widetilde{F}_{T_{1}}}\lesssim(\|q_{1}\|_{\widetilde{L}^{2}_{T_{1}}(B^{\NN+1}_{p,1})}+\|q_{2}\|_{\widetilde{L}_{T_{1}}^{2}(B^{\NN+1}_{p,1})}+
\|u_{1}\|_{\widetilde{L}^{2}_{T_{1}}(B^{\NN}_{p,1})}+\|u_{2}\|_{\widetilde{L}_{T_{1}}^{2}(B^{\NN}_{p,1})})\\
&\hspace{11cm}\times\|(\de
q,\de v)\|_{\widetilde{F}_{T_{1}}}.
\end{aligned}
$$
We have then for $T_{1}$ small enough: $(\de q,\de v)=(0,0)$ on $[0,T_{1}]$ and by connectivity we finally
conclude that:
$$q_{1}=q_{2},\;v_{1}=v_{2}\;\;\mbox{on}\;\;[0,T].$$ \hfill {$\Box$}
\subsection{Estimate of the time of existence}

Here we can estimate $T$ such that the assumption $({\cal{H}}_{\e})$ is verified, indeed we have seen in proposition \ref{chaleur} that when $1\leq\rho_1\leq+\infty$ it exists two constants $c,C>0$ such that $v_L$ verifying a heat equation has the following property:
\begin{equation}
\|v_L\|_{\widetilde{L}^{\rho_1}_{T}(B^{s+\frac{2}{\rho_1}}_{p,1})}\leq C\big(\sum_{q\in\mathbb{Z}}2^{qs}\|\D_q v_0\|_{L^p}\big(\frac{1-e^{-c\mu T2^{2q}\rho_1}}{c\mu\rho_1})^{\frac{1}{\rho_1}}\big).
\label{estimtemps}
\end{equation}
In particular we have:
\begin{equation}
\|v_L\|_{\widetilde{L}^{1}_{T}(B^{\NN+1}_{p,1})}\leq C\big(\sum_{q\in\mathbb{Z}}2^{q(\NN-1)}\|\D_q v_0\|_{L^p}\big(\frac{1-e^{-c\mu T2^{2q}}}{c\mu})\big).
\label{estimtemps1}
\end{equation}
Now estimating $q_L$ in the system (\ref{lineaire}) and using again the proposition (\ref{chaleur}), we get for $C>0$ large enough :
\begin{equation}
\begin{aligned}
&\|q_L\|_{\widetilde{L}^{1}_{T}(B^{\NN+2}_{p,1})}\leq C\big(\sum_{q\in\mathbb{Z}}2^{q\NN}\|\D_q q_0\|_{L^p}\big(\frac{1-e^{-c\mu T2^{2q}}}{c\mu})\\
&\hspace{2cm}+\sum_{q\in\mathbb{Z}}2^{q\NN}\|\D_q {\rm div}v_L\|_{L^1_T(L^p)}\big(\frac{1-e^{-c\mu T2^{2q}}}{c\mu})\big),\\
&\leq  C\big(\sum_{q\in\mathbb{Z}}2^{q\NN}\|\D_q q_0\|_{L^p}\big(\frac{1-e^{-c\mu T2^{2q}}}{c\mu})+
\|{\rm div}v_L\|_{L^1_T(B^{\NN}_{p,1})} \big)
\end{aligned}
\label{estimtemps2}
\end{equation}
We deduce that for $C$ large enough:
\begin{equation}
\begin{aligned}
&\|q_L\|_{\widetilde{L}^{1}_{T}(B^{\NN+2}_{p,1})}\leq C\big(\sum_{q\in\mathbb{Z}}2^{q\NN}\|\D_q q_0\|_{L^p}\big(\frac{1-e^{-c\mu T2^{2q}}}{c\mu})\\
&\hspace{2cm}+ \sum_{q\in\mathbb{Z}}2^{q(\NN-1)}\|\D_q v_0\|_{L^p}\big(\frac{1-e^{-c\mu T2^{2q}}}{c\mu})  \big)
\end{aligned}
\label{estimtemps3}
\end{equation}
It remains only to choose $T$ sufficiently small such that $( {\cal{H}}_{\e})$ is verified, let us start with the case:
$$ C\sum_{q\in\mathbb{Z}}2^{q\NN}\|\D_q q_0\|_{L^p}\big(\frac{1-e^{-c\mu T2^{2q}}}{c\mu})\leq \frac{\e}{2}.$$
We have then using the fact that $1-e^{-x}\leq x$ when $x\in\R$:
$$
\begin{aligned}
&\sum_{q\in\mathbb{Z}}2^{q\NN}\|\D_q q_0\|_{L^p}\big(\frac{1-e^{-c\mu T2^{2q}}}{c\mu})\leq \sum_{q\leq l_0}2^{q\NN}\|\D_q q_0\|_{L^p} T2^{2q}+\frac{2}{c\mu}\sum_{q\geq l_0}2^{q\NN}\|\D_q q_0\|_{L^p},\\
&\leq (2^{(2-\e')l_0}T+\frac{2}{c\mu} 2^{-l_0 \e'})\|q_0\|_{B^{\NN+\e'}_{p,1}}.
\end{aligned}
$$
Let us choose $l_0$ such that:
\begin{equation}
C\frac{2}{c\mu} 2^{-l_0 \e'}\|q_0\|_{B^{\NN+\e'}_{p,1}}\leq\frac{\e}{4}.
\end{equation}
In particular we must have:
\begin{equation}
2^{l_0 \e'}\geq\frac{8C \|q_0\|_{B^{\NN+\e'}_{p,1}}}{c\mu \e}.
\end{equation}
In particular it implies that:
$$l_0\geq \frac{\frac{1}{\e'}\ln(\frac{8C \|q_0\|_{B^{\NN+\e'}_{p,1}}}{c\mu\e})}{\ln 2}.$$
Let us choose now $l_0= \frac{\frac{1}{\e'}\ln(\frac{8C \|q_0\|_{B^{\NN+\e'}_{p,1}}}{c\mu\e})}{\ln 2}$ and we want to ensure that:
$$C2^{(2-\e')l_0}T \|q_0\|_{B^{\NN+\e'}_{p,1}}\leq \frac{\e}{4}.$$
In particular we must have:
$$
\begin{aligned}
&T\leq\frac{\e}{4C \|q_0\|_{B^{\NN+\e'}_{p,1}}2^{(2-\e')l_0}}.
\end{aligned}
$$
Since we have:
$$2^{(2-\e')l_0}=\big(\frac{8C \|q_0\|_{B^{\NN+\e'}_{p,1}}}{c\mu\e}\big)^{\frac{2}{\e'}} \frac{c\mu\e}{8C \|q_0\|_{B^{\NN+\e'}_{p,1}}} $$
It implies that:
$$
\begin{aligned}
&T\leq \frac{ 2(c\mu)^{\frac{2}{\e'}-1}   \e^{\frac{2}{\e'}}}{(8C)^{\frac{2}{\e'}}\|q_0\|_{B^{\NN+\e'}_{p,1}}^{\frac{2}{\e'}}}
\end{aligned}
$$
In other term we have prove that for:
\begin{equation}
T=\min \big(\frac{ 2(c\mu)^{\frac{2}{\e'}-1}   \e^{\frac{2}{\e'}}}{(8C)^{\frac{2}{\e'}}\|q_0\|_{B^{\NN+\e'}_{p,1}}^{\frac{2}{\e'}}},\frac{ 2(c\mu)^{\frac{2}{\e'}-1}   \e^{\frac{2}{\e'}}}{(8C)^{\frac{2}{\e'}}\|v_0\|_{B^{\NN-1+\e'}_{p,1}}^{\frac{2}{\e'}}},
\big),
\label{3.34}
\end{equation}
the assumption $({\cal H}_{\e})$ is verified. We recall that we need also that $T$ verifies the condition (\ref{impborneT}), so we have finally:
\begin{equation}
\begin{aligned}
&T=D(\|q_0\|_{B^{\NN}_{p,1}}, \|v_0\|_{B^{\NN-1}_{p,1}},\|q_0\|_{B^{\NN+\e'}_{p,1}}, \|v_0\|_{B^{\NN-1+\e'}_{p,1}}  )\\
&\leq   \big(\frac{ 2(c\mu)^{\frac{2}{\e'}-1}   \e^{\frac{2}{\e'}}}{(8C)^{\frac{2}{\e'}}\|q_0\|_{B^{\NN+\e'}_{p,1}}^{\frac{2}{\e'}}},\frac{ 2(c\mu)^{\frac{2}{\e'}-1}   \e^{\frac{2}{\e'}}}{(8C)^{\frac{2}{\e'}}\|v_0\|_{B^{\NN-1+\e'}_{p,1}}^{\frac{2}{\e'}}},\\
&\hspace{2cm}\frac{1}{16 C^2_1  (\|q_0\|_{B^{\NN}_{p,1}}+\|v_0\|_{B^{\NN-1}_{p,1}})(1+\sqrt{ \|q_0\|_{B^{\NN}_{p,1}}+\|v_0\|_{B^{\NN-1}_{p,1}}})^2)},\frac{C_1}{4}).
\end{aligned}
\label{3.34}
\end{equation}
Here $D$ is a decreasing function in terms of each variable. It implies in particular that the time of existence $T^{*}$ of the strong solution is superior at least to $T$.
\subsection{Propagation of the regularity}
\label{regu}
Here we assume in addition that $(q_0,v_0)$ belongs in $B^{s}_{p,1}\times B^{s-1}_{p,1}$ with $s>\NN$. In a very classical way the regularity is preserved on $[0,T]$, it means that:
$$(q,u)\in \widetilde{C}_{T}(B^{s}_{p,1})\times  \widetilde{C}_{T}(B^{s-1}_{p,1}).$$
It suffices to proceed as in the previous section.
\section{Uniform energy estimates in the particular case $P(\rho)=a\rho$}
\label{section4}
In the next sections, we are going to prove the theorem \ref{theo2} in the particular case of $P(\rho)=a\rho$. In the section \ref{section6}, we will give a proof for general pressure following the assumptions of the theorem \ref{theo2}.\\
In this section we assume that $(q_0,v_0)$ verify the assumption of the theorem \ref{theo2}. In particular since $(q_0,v_0)$ belongs in $B^{\frac{N}{p}}_{p,1}\times B^{\frac{N}{p}-1}_{p,1}$ for $N<p<2N$, we know via the theorem \ref{theo1} that it exists a strong solution $(q,v)$ to the system (\ref{systK1}) on the interval $(0,T^{*})$ with $T^*$ the maximal time of existence. In the sequel we shall assume by absurd that $T^*<+\infty$, our goal is now to prove that it is not possible.\\
In addition we assume an extra hypothesis on the initial data since we choose $(q_0,v_0)$ in $B^s_{p,1}\times B^{s-1}_{p,1}$ with $s>\NN$ sufficiently large such that our solution $(q,v)$ are in $C^k((0,T^*)\times\R^N)$ in order to justify all the integrations by parts. Let us mention that this regularity result is obtain by Besov embedding and the propagation of the regularity proved in theorem \ref{theo1}. More precisely we have for any $T\in(0,T^*)$ via the theorem \ref{theo1}:
\begin{equation}
\begin{cases}
\begin{aligned}
&q\in\widetilde{C}([0,T], B^{\NN}_{p,1}\cap B^{s}_{p,1}  )\cap L^1([0,T],B^{\NN+2}_{p,1}\cap B^{s+2}_{p,1}),\\
&v\in\widetilde{C}([0,T], B^{\NN-1}_{p,1}\cap B^{s-1}_{p,1}  )\cap L^1([0,T],B^{\NN+1}_{p,1}\cap B^{s+1}_{p,1}),\\
&\rho,\frac{1}{\rho}\in L^{\infty}_T(L^{\infty}(\R^N)).
\end{aligned}
\end{cases}
\label{reg}
\end{equation}
In addition in order to justify that $v$ belongs in $C([0,T], L^{p_1})$ for any $T\in(0,T^*)$
we assume in addition that $(q_0,v_0)$ are in all $B^{1}_{p_1,1}\times B^0_{p_1,1}$ with $2\leq p_1<+\infty$.
\begin{remarka}
In the theorem \ref{theo2} the assumptions on the initial data are more restrictive, it suffices then when we have prove the result of global strong solution for the previous choice of initial data to use a regularization process and to pass to the limit using compactness argument.
\label{remtech}
\end{remarka}
We are now going to prove a series of estimates on $\rho$ and $v$ which will depend only on the time $T$. All the integration by parts and all the quantity estimates
 will have sense via the condition (\ref{reg}) and the remarks \ref{remtech}.
%
\subsection{Energy estimates}
Following Lions in \cite{fL2} p 207 we set: $q(t)=t\int^{t}_{1}\frac{P(s)}{s^2}ds$. Since $P(\rho)=a\rho$ we have $q(t)=at\ln t$. Multiplying the momentum equation of (\ref{systK}) by $v$ we get:
\begin{equation}
\begin{aligned}
&\int_{\R^N}\big(\frac{1}{2}\rho(t,x)|v(t,x)|^{2}+(\Pi(\rho)(t,x)-\Pi(\bar{\rho}))\big)dx+\int^{t}_{0}\int_{\R^N}\big(\mu\,\rho(t,x)|\n v|^{2}(t,x)\\
&\hspace{2cm}+\frac{a\mu}{\rho}|\n\rho|^{2}(t,x)\big)dtdx \leq \int_{\R^N}\big(\rho_{0}(x)|v_{0}(x)|^{2}+(\Pi(\rho_{0})(x)-\Pi(\bar{\rho}))\big)dx.
\end{aligned}
\label{energie}
\end{equation}
with $\Pi(\rho)$ defined as follows:
$$\Pi(\rho)=a\big(\rho\ln(\frac{\rho}{\bar{\rho}})+\bar{\rho}-\rho\big)=q(\rho)-q(\bar{\rho})-q'(\bar{\rho})(\rho-\bar{\rho}).$$
We have then $\Pi(\bar{\rho})=\Pi'(\bar{\rho})=0$ with:
$$\Pi'(s)=a\ln(\frac{s}{\bar{\rho}}),\;\;\Pi''(s)=\frac{a}{s}.$$
We deduce that $\Pi$ is convex. It implies in particular that $(\Pi(\rho)-\Pi(\bar{\rho}))\geq 0$ and using
(\ref{energie}) we observe that $(\Pi(\rho)-\Pi(\bar{\rho}))$ is in $L^\infty_T(L^1(\R^N))$ for any $T>0$.\\
Let us prove now useful proposition which ensure a $H^1$ control on the density. 
\begin{proposition}
We have $(\rho-\bar{\rho})\in L^\infty(L^{2}_{1}(\R^N))$
\end{proposition}
{\bf Proof:} For $\delta>0$ we show that it exists $C>0$ such that:
$$\frac{1}{C}|\rho-\bar{\rho}|1_{\{ |\rho-\bar{\rho}|\geq \delta\} }\leq (\Pi(\rho)-\Pi(\bar{\rho}))1_{\{ |\rho-\bar{\rho}|\geq \delta\} }.$$
Next since $\Pi(\bar{\rho})=\Pi'(\bar{\rho})=0$ it exists $C>0$ such that:
$$\frac{1}{C}|\rho-\bar{\rho}|^21_{\{ |\rho-\bar{\rho}|\leq \delta\} }\leq (\Pi(\rho)-\Pi(\bar{\rho}))1_{\{ |\rho-\bar{\rho}|\leq \delta\} }.$$
It implies in particular that $(\rho-\bar{\rho})$ is in $L^\infty_T(L^1(\R^N))$ for any $T>0$ using the fact that $(\Pi(\rho)-\Pi(\bar{\rho}))$ is in $L^\infty_T(L^1(\R^N))$. $\blacksquare$
\begin{proposition}
We have $(\sqrt{\rho}-\sqrt{\bar{\rho}})\in L^\infty(L^2).$
\end{proposition}
{\bf Proof:} We have:
$$(\sqrt{\rho}-\sqrt{\bar{\rho}})^2=(\rho-\bar{\rho})+2\sqrt{\bar{\rho}}(\sqrt{\bar{\rho}}-\sqrt{\rho}).$$
By Young inequality we have:
$$\frac{1}{2}|\sqrt{\rho}-\sqrt{\bar{\rho}}|^2\leq |\rho-\bar{\rho}|+2\bar{\rho}$$
And we have then:
$$\frac{1}{2}\|(\sqrt{\rho}-\sqrt{\bar{\rho}})1_{\{|\rho-\bar{\rho}|\geq\delta\}}\|^2_{L^2}\leq \|(\rho-\bar{\rho})   1_{\{|\rho-\bar{\rho}  |\geq\delta\}}\|_{L^1}+2\bar{\rho}|  \{|\rho-\bar{\rho}  |\geq\delta\}|<+\infty,$$
because $(\rho-\bar{\rho})\in L^\infty(L^2_1(\R^N))$. Next it exists $C>0$ such that:
$$\frac{1}{C}|\sqrt{\rho}-\sqrt{\bar{\rho}}|1_{\{|\rho-\bar{\rho}|\leq\delta\}}\leq 
\frac{1}{2\sqrt{\bar{\rho}}}|\rho-\bar{\rho}|1_{\{|\rho-\bar{\rho}|\leq\delta\}} .$$
And we deduce that $(\sqrt{\rho}-\sqrt{\bar{\rho}})1_{\{|\rho-\bar{\rho}|\leq\delta\}}$ is in $L_T^\infty(L^2)$ since $(\rho-\bar{\rho})$ is in $L^\infty_T(L^2_1(\R^N))$. It concludes the proof of the proposition. $\blacksquare$ 
\begin{proposition}
We have $(\sqrt{\rho}-\sqrt{\bar{\rho}})\in L^\infty(H^1(\R^N))$.
\label{22energy}
\end{proposition}
{\bf Proof:} Since we have shown that $(\sqrt{\rho}-\sqrt{\bar{\rho}})$ is in $L^\infty_T(L^2(\R^N))$ and that $\n\sqrt{\rho}$ is also in $L^\infty_T(L^2(\R^N))$ via the energy estimate (\ref{energie}) it concludes the proof of the proposition. $\blacksquare$ 
\subsection{Gain of integrability on $v$}
In this part, we are going to prove that the effective velocity $v$ preserves in some sense the $L^p$ norm following an idea developed by Mellet and Vasseur in \cite{fMV1,5MV2} for the compressible Navier-Stokes system. More precisely we are going to show that $\rho^{\frac{1}{p}} v$ is in $L^{\infty}((0,T),L^{p}(\R^{N}))$
for any $2\leq p<+\infty$ provided that $\rho^{\frac{1}{p}} v_0$ is also in $L^p(\R^{N})$ for any $2\leq p<+\infty$. We shall observe in the proof of this result that it is strongly related to our choice of pressure $P(\rho)=a\rho$ with $a>0$. 
\begin{lemme}
\label{gainin}
Let $(\rho,v)$ be  our strong solution  on $(0,T^*)$, then it exists $C$ an increasing function  depending only on the initial data of theorem \ref{theo2}  such that for all $T\in(0,T^*)$ we have for all $p\in[4,+\infty)$:
$$
\begin{aligned}
&\|\rho^{\frac{1}{p}}v(T,\cdot)\|_{L^p}\leq C(T).
\end{aligned}
$$
\end{lemme}
\begin{remarka}
Let us point out that $C(T)$ does not depend on $p\in[4,+\infty)$.
\end{remarka}
{\bf Proof:} As in \cite{fMV1,5MV2}, we now want to obtain additional information on the integrability of $v$, to do it  we multiply the momentum equation of (\ref{systK}) by $v|v|^{p-2}$ and integrate over $\R^{N}$, we obtain then:
\begin{equation}
\begin{aligned}
&\frac{1}{p}\int_{\R^{N}}\rho\p_{t}(|v|^{p})dx+\int_{\R^{N}}\rho u\cdot\n(\frac{|v|^{p}}{p})dx
+\int_{\R^{N}} \rho|v|^{p-2}|\n v|^{2}dx\\
&\hspace{1,5cm}+(p-2)\int_{\R^{N}} \rho\sum_{i,j,k}v_{j}v_{k}\p_{i}v_{j}\p_{i}v_{k}|v|^{p-4}dx+\int_{\R^{N}} |v|^{p-2}v\cdot\n (a\rho )dx=0.
\end{aligned}
\label{in1A}
\end{equation}
Next we observe that:
$$\sum_{i,j,k}v_{j}v_{k}\p_{i}v_{j}\p_{i}v_{k}=\sum_{i}(\sum_{j}v_{j}\p_{i}v_{j})^{2}=\sum_{i}\big[\frac{1}{2}\p_{i}(|v|^{2})\big]^{2}=\frac{1}{4}|\n(|v|^{2})|^{2}.$$
We get then as ${\rm div}(\rho u)=-\p_{t}\rho$ and by using (\ref{in1A}):
\begin{equation}
\begin{aligned}
&\frac{1}{p}\int_{\R^{N}}\p_{t}(\rho|v|^{p})dx+\int_{\T^{N}} \rho|v|^{p-2}|\n v|^{2}dx+\frac{(p-2)}{4}\int_{\R^{N}} \rho|\n(|v|^{2})|^{2}|v|^{p-4}dx\\
&\hspace{8cm}+\int_{\R^{N}} |v|^{p-2}v\cdot\n (a\rho) dx=0.
\end{aligned}
\label{in1A1}
\end{equation}
We have then by integrating over $(0,t)$ with $0<t\leq T$:
\begin{equation}
\begin{aligned}
&\frac{1}{p}\int_{\R^{N}}(\rho|v|^{p})(t,x)dx+\int^{t}_{0}\int_{\R^{N}} \rho|v|^{p-2}|\n v|^{2}(s,x)dsdx\\
&+\frac{(p-2)}{4}\int^{t}_{0} \int_{\R^{N}}\rho|\n(|v|^{2})|^{2}|v|^{p-4}(s,x)dsdx\leq \frac{1}{p}\int_{\R^{N}}(\rho_{0}|v_{0}|^{p})(x)dx\\
&\hspace{7cm}+|\int^{t}_{0} \int_{\R^{N}}|v|^{p-2}v\cdot\n a\rho(s,x)dsdx|.
\end{aligned}
\label{in1A1}
\end{equation}
By integration by parts we have since $(\rho,v)$ are regular:
$$\int^{t}_{0} \int_{\R^{N}}|v|^{p-2}v\cdot\n a\rho(s,x)dsdx=-a\int^{t}_{0} \int_{\R^{N}}{\rm div}(|v|^{p-2}v)\rho(s,x)dsdx.$$
Next we have:
$${\rm div}(|v|^{p-2}v)=|v|^{p-2}{\rm div}(v)+(p-2)|v|^{p-4}v\cdot(v\cdot\n v).$$
By Young inequality we have:
\begin{equation}
\begin{aligned}
&|-a\int^{t}_{0} \int_{\R^{N}}|v|^{p-2}{\rm div}(v)\rho(s,x)dsdx|\lesssim \frac{\e}{2}\int^{t}_{0} \int_{\R^{N}}\rho|v|^{p-2}|\n v|^2 dsdx\\
&\hspace{7cm}+\frac{1}{2\e}a^2N^2\int^{t}_{0} \int_{\R^{N}}\rho
|v|^{p-2}dsdx.
\end{aligned}
\label{Y1}
\end{equation}
Plugging (\ref{Y1}) in (\ref{in1A1}) with $\e=1$ we have:
\begin{equation}
\begin{aligned}
&\frac{1}{p}\int_{\R^{N}}(\rho|v|^{p})(t,x)dx+\frac{1}{2}\int^{t}_{0}\int_{\R^{N}} \rho|v|^{p-2}|\n v|^{2}(s,x)dsdx\\
&+\frac{(p-2)}{4}\int^{t}_{0} \int_{\R^{N}}\rho|\n(|v|^{2})|^{2}|v|^{p-4}(s,x)dsdx\leq \frac{1}{p}\int_{\R^{N}}(\rho_{0}|v_{0}|^{p})(x)dx\\
&\hspace{2cm}+\frac{1}{2}a^2N^2\int^{t}_{0} \int_{\R^{N}}\rho
|v|^{p-2}dsdx+(p-2)|\int^t_0\int_{\R^N}\rho |v|^{p-4}v\cdot(v\cdot\n v) dx ds|.
\end{aligned}
\label{1in1A1}
\end{equation}
It remains to estimates the two last terms on the right hand side of (\ref{1in1A1}). We have for $p\geq 4$:
$$\rho|v|^{p-2}=|\sqrt{\rho}v|^{2(1-\theta)}|\rho^{\frac{1}{p}}v|^{p\theta},$$
with $p-2=2(1-\theta)+p\theta$ with $\theta\in(0,1)$. We deduce by H\"older's inequality that:
$$  \int_{\R^{N}}\rho|v|^{p-2}dx\leq\|\sqrt{\rho}v\|_{L^2}^{2(1-\theta)}\|\rho^{\frac{1}{p}}v\|_{L^p}^{p\theta}. $$
In particular we have:
$$\theta=\frac{p-4}{p-2},\;1-\theta=\frac{2}{p-2}. $$
Then by H\"older's inequality we get:
$$  
\begin{aligned}
&\int^t_0 \int_{\R^{N}}\rho|v|^{p-2}dx\leq \|\sqrt{\rho}v \|_{L_t^\infty(L^2)}^{\frac{4}{p-2}} \int^t_0 \|\rho^{\frac{1}{p}}v(s)\|_{L^p}^{\frac{p(p-4)}{p-2}}ds.
\end{aligned}
$$
We deduce by Young inequality (and choosing $\e=\frac{1}{p}$) that:
\begin{equation}  
\begin{aligned}
&\int^t_0 \int_{\R^{N}}\rho|v|^{p-2}dx\leq \|\sqrt{\rho}v \|_{L_t^\infty(L^2)}^{\frac{4}{p-2}} \int^t_0 \|\rho^{\frac{1}{p}}v\|_{L^p}^{\frac{p(p-4)}{p-2}}ds,\\
&\leq \|\sqrt{\rho}v \|_{L_t^\infty(L^2)}^{\frac{4}{p-2}}  \int^t_0\big(\frac{\e(p-4)}{p-2} \|\rho^{\frac{1}{p}}v\|_{L^p}^{p}+\frac{2}{\e (p-2)}\big)ds,\\
&\leq \|\sqrt{\rho}v \|_{L_t^\infty(L^2)}^{\frac{4}{p-2}}  \big(\frac{(p-4)}{p-2} \int^t_0\frac{1}{p}\|\rho^{\frac{1}{p}}v(s,\cdot)\|_{L^p}^{p} ds+\frac{2p}{ (p-2)}t\big).
\end{aligned}
\label{etim1}
\end{equation}
Similarly we have by Young inequality:
\begin{equation}
\begin{aligned}
&|\int^t_0\int_{\R^N}(p-2)a\rho |v|^{p-4}v\cdot(v\cdot\n v)dxds|= |\int^t_0\int_{\R^N}\frac{(p-2)}{2}a\rho |v|^{p-4}v\cdot\n|v|^2 dxds|,\\
&\lesssim\frac{\e}{2} \int^t_0\int_{\R^N}\rho |v|^{p-4}|\n|v|^2|^2dxds+\frac{1}{2\e}\frac{(p-2)^2 a^2}{4}\int^t_0\int_{\R^N}\rho |v|^{p-2}dx\,ds.
\end{aligned}
\label{Y2}
\end{equation}
We choose $\e=\frac{p-2}{4}$ and plugging (\ref{Y2}) in (\ref{1in1A1}) we have:
\begin{equation}
\begin{aligned}
&\frac{1}{p}\int_{\R^{N}}(\rho|v|^{p})(t,x)dx+\frac{1}{2}\int^{t}_{0}\int_{\R^{N}} \rho|v|^{p-2}|\n v|^{2}(s,x)dsdx\\
&+\frac{(p-2)}{8}\int^{t}_{0} \int_{\R^{N}}\rho|\n(|v|^{2})|^{2}|v|^{p-4}(s,x)dsdx\leq \frac{1}{p}\int_{\R^{N}}(\rho_{0}|v_{0}|^{p})(x)dx\\
&\hspace{4cm}+[\frac{1}{2}a^2N^2+\frac{a^2}{2}(p-2)]\int^{t}_{0} \int_{\R^{N}}\rho
|v|^{p-2}dsdx.
\end{aligned}
\label{1in1A2}
\end{equation}
Next using (\ref{etim1}) we have with $\e=\frac{1}{p(p-4)}$:
\begin{equation}
\begin{aligned}
&(p-2)\int^t_0 \int_{\R^{N}}\rho|v|^{p-2}dxds\\
&\leq \|\sqrt{\rho}v \|_{L_t^\infty(L^2)}^{\frac{4}{p-2}}  \big(\e(p-4) \int^t_0\|\rho^{\frac{1}{p}}v\|_{L^p}^{p} ds+\frac{2}{\e}t\big)\\
&\leq \|\sqrt{\rho}v \|_{L_t^\infty(L^2)}^{\frac{4}{p-2}}  \big(\frac{1}{p} \int^t_0\|\rho^{\frac{1}{p}}v\|_{L^p}^{p} ds+2p (p-4)t\big)
\end{aligned}
\label{etim2}
\end{equation}
Finally plugging (\ref{etim1}) and (\ref{etim2}) in (\ref{1in1A2}) we have:
\begin{equation}
\begin{aligned}
&\frac{1}{p}\int_{\R^{N}}(\rho|v|^{p})(t,x)dx+\frac{1}{2}\int^{t}_{0}\int_{\R^{N}} \rho|v|^{p-2}|\n v|^{2}(s,x)dsdx\\
&+\frac{(p-2)}{8}\int^{t}_{0} \int_{\R^{N}}\rho|\n(|v|^{2})|^{2}|v|^{p-4}(s,x)dsdx\leq \frac{1}{p}\int_{\R^{N}}(\rho_{0}|v_{0}|^{p})(x)dx\\
&+\|\sqrt{\rho}v \|_{L_t^\infty(L^2)}^{\frac{4}{p-2}}  \biggl(\frac{a^2}{2}(N^2\frac{(p-4)}{p-2}+1) \int^t_0\frac{1}{p}\|\rho^{\frac{1}{p}}v(s,\cdot)\|_{L^p}^{p} ds+ \frac{a^2 t}{2}\big(N^2\frac{2p}{ (p-2)}+2p (p-4)\big)\biggl).
\end{aligned}
\label{in1A1a}
\end{equation}
By Gronwall lemma, we conclude that for all $t\in[0,T]$ we have:
$$
\begin{aligned}
&\frac{1}{p}\int_{\R^{N}}(\rho|v|^{p})(t,x)dx\leq  \biggl(\frac{1}{p}\int_{\R^{N}}(\rho_{0}|v_{0}|^{p})(x)dx+\|\sqrt{\rho}v \|_{L_T^\infty(L^2)}^{\frac{4}{p-2}}\frac{a^2}{2}\big(N^2 \frac{2p}{ (p-2)}+2p (p-4)\big)T\biggl)\\
&\hspace{4cm}\times \exp(\|\sqrt{\rho}v \|_{L_T^\infty(L^2)}^{\frac{4}{p-2}}\frac{a^2}{2}(N^2\frac{(p-4)}{p-2}+1)t).
\end{aligned}
$$
In particular we have:
\begin{equation}
\begin{aligned}
&\|\rho^{\frac{1}{p}}v(t,\cdot)\|_{L^p}\leq  \biggl(\int_{\R^{N}}(\rho_{0}|v_{0}|^{p})(x)dx+\|\sqrt{\rho}v \|_{L_T^\infty(L^2)}^{\frac{4}{p-2}}\frac{a^2}{2}\big( N^2 \frac{2p^2}{ (p-2)}+2p^2 (p-4)\big)T\biggl)^{\frac{1}{p}}\\
&\hspace{4cm}\times \exp(\frac{1}{p}\|\sqrt{\rho}v \|_{L_T^\infty(L^2)}^{\frac{4}{p-2}}\frac{a^2}{2}(N^2 \frac{(p-4)}{p-2}+1)t).
\end{aligned}
\label{cruciala}
\end{equation}
For $x,y \geq 0$ and for all $p\geq 1$ we verify easily that:
$$(x+y)^{\frac{1}{p}}\leq 2^{\frac{1}{p}}(x^{\frac{1}{p}}+y^{\frac{1}{p}})$$
then for all $p\geq 4$ we have:
\begin{equation}
\begin{aligned}
&\|\rho^{\frac{1}{p}}v(t,\cdot)\|_{L^p}\leq 2^{\frac{1}{p}}\biggl(\|\rho_0^{\frac{1}{p}}v_0\|_{L^p}+\|\sqrt{\rho}v \|_{L_T^\infty(L^2)}^{\frac{4}{p(p-2)}}(\frac{a^2}{2})^{\frac{1}{p}} \big(N^2 \frac{2p^2}{ (p-2)}+2p^2 (p-4)\big)^{\frac{1}{p}}T^{\frac{1}{p}}\biggl)\\
&\hspace{4cm}\times \exp(\frac{1}{p}\|\sqrt{\rho}v \|_{L_T^\infty(L^2)}^{\frac{4}{p-2}}\frac{a^2}{2}(N^2\frac{(p-4)}{p-2}+1)t).
\end{aligned}
\label{crucialb}
\end{equation}
It ends up the proof of the lemma \ref{gainin}.
$\blacksquare$
\subsection{Control on $\max_{(t,x)\in [0,T]\times \R^N}\frac{1}{\rho (t,x)}$ for $0<T<T^*$}
\begin{proposition}
Under the assumption of theorem \ref{theo2}, the density $\rho$ verifies for any $T\in (0,T^*)$ with in our case $T^*<+\infty$ which is a priori strictly finite:
\begin{equation}
\|\frac{1}{\rho}\|_{L^\infty_T(L^\infty(\R^N))}\leq C(T),
\label{contvide}
\end{equation}
with $C$ an increasing function in $T$.
\label{avide}
\end{proposition}
{\bf Proof}Let us deal with the momentum equation:
$$
\begin{cases}
\begin{aligned}
&\p_t\rho-\mu\D\rho+{\rm div}(\rho v)=0,\\
&\rho(0,\cdot)=\rho_0.
\end{aligned}
\end{cases}
$$
On $(0,T^*)$ since the density $\rho$ is regular and is far away from the vacuum, $\rho$ verifies the following equation:
$$\p_t(\frac{1}{\rho^{\alpha}})-\alpha\rho^{-\alpha-1}{\rm div}(\rho u)=0.$$
And we deduce since $v=u+\mu\n\ln\rho$ that:
$$\p_t(\rho^{-\alpha})-\alpha\rho^{-\alpha}{\rm div}v+v\cdot\n\rho^{-\alpha}+\alpha\mu\rho^{-\alpha-1}\D\rho =0.$$
Since
$$-\alpha \rho^{-\alpha-1}\D\rho=\D\rho^{-\alpha}-\alpha(\alpha+1)\rho^{-\alpha-2}|\n\rho|^2,$$
we obtain then the following equation:
\begin{equation}
\p_t\rho^{-\alpha}-\mu\D\rho^{-\alpha}+\mu\alpha(\alpha+1)\rho^{-\alpha-2}|\n\rho|^2+v\cdot\n\rho^{-\alpha}-\alpha\rho^{-\alpha}{\rm div}v=0.
\label{systcru1}
\end{equation}
We are going following in a crucial way some deep ideas due to Ladyzenskaja et al in \cite{La} which extend the De Giorgi method (introduced to obtain regularity results for the elliptic equations) to the parabolic equations. Let us start by multiplying the previous equation  by $(\rho^{-\alpha})^{(k)}=\max (\rho^{-\alpha}(t,x)-k,0),$ and integrate over $(0,t_1)\times\R^N$, it gives using the fact that the term $\mu\alpha(\alpha+1)\rho^{-\alpha-2}|\n\rho|^2$ is positive: 
\begin{equation}
\begin{aligned}
&\frac{1}{2}\int_{\R^N} [(\rho^{-\alpha})^{(k)}(t_1,x)]^2 dx+\mu\int^{t_1}_0\int_{A_k(t)}|\n\rho^{-\alpha} (t,x)|^2 dx dt\\
&\hspace{0,5cm}\leq -(\alpha+1) \int^{t_1}_0\int_{A_k(t)} v_i \p_i (\rho^{-\alpha}) (\rho^{-\alpha})^{(k)}(t,x) dx dt -\alpha\int^{t_1}_0\int_{A_k(t)}\rho^{-\alpha} v_i \p_i (\rho^{-\alpha}) dx dt.
\end{aligned}
\label{impcru}
\end{equation}
with:
$$A_k(t)=\{ x\in\R^N;\;\frac{1}{\rho^{\alpha}(t,x)}\geq k\}=\{ x\in\R^N;\;0\leq\rho(t,x)\leq\frac{1}{k^{\frac{1}{\alpha}}}\}.$$
\begin{remarka}
Let us mention that $\rho$ is regular choosing initial data $(\rho_0-\bar{\rho})$ sufficiently regular in Besov space and using the fact that the regularity is preserved on $(0,T^*)$ by the section \ref{regu}. In addition the integrals in (\ref{impcru}) are well defined, indeed if we consider:
$$\int_{\R^N} [(\rho^{-\alpha})^{(k)}(t_1,x)]^2 dx,$$
we know that $\mbox{supp}\big[(\rho^{-\alpha})^{(k)}\big]$  is included in $A_k(t_1)$. It implies that  $\mbox{supp}(\rho^{-\alpha})^{(k)}$ is in $\{|\rho(t_1,x)-\bar{\rho}|\geq (\bar{\rho}-\frac{1}{k^{\frac{1}{\alpha}}})\}$. This set is of finite measure since we know that $\rho(t_1,\cdot)$ belongs in $L^2_1(\R^N)$, and using the fact that for all $t\in (0,T^*)$ $\frac{1}{\rho(t,\cdot)}$ is in $L^\infty(\R^N)$ we deduce that the previous integral is finite.
\end{remarka}
Now by Young inequality we have:
\begin{equation}
\begin{aligned}
&\frac{1}{2}\int_{\R^N} [(\rho^{-\alpha})^{(k)}(t_1,x)]^2 dx+\mu\int^{t_1}_0\int_{A_k(t)}|\n\rho^{-\alpha} (t,x)|^2 dx dt\\
&\leq |(\alpha+1)| \int^{t_1}_0\int_{A_k(t)}\big( \sum_i \frac{1}{2\e}| v_i|^2 (\rho^{-\alpha}-k)^2 +\frac{\e}{2}| \n (\rho^{-\alpha})|^2\big) (t,x) dx dt \\
&+|\alpha|\int^{t_1}_0\int_{A_k(t)}\big(\frac{1}{2\e_1}|\rho^{-\alpha}|^2\sum_i |v_i|^2 +\frac{\e_1}{2}|\n  (\rho^{-\alpha})|^2\big)(t,x) dx dt.
\end{aligned}
\label{impcru1}
\end{equation}
We choose $\e=\frac{\mu}{4|\alpha+1|}$, $\e_1=\frac{\mu}{4|\alpha|}$ and we have for $k\geq 1$:
\begin{equation}
\begin{aligned}
&\frac{1}{2}\int_{\R^N} [(\rho^{-\alpha})^{(k)}(t_1,x)]^2 dx+\frac{\mu}{2}\int^{t_1}_0\int_{A_k(t)}|\n\rho^{-\alpha} (t,x)|^2 dx dt\\
&\leq \frac{2|(\alpha+1)|^2}{\mu} \int^{t_1}_0\int_{A_k(t)}\big( |v|^2 [(\rho^{-\alpha}-k)^2+k^2] \big) (t,x) dx dt \\
&+\frac{2|\alpha|^2}{\mu}\int^{t_1}_0\int_{A_k(t)}\big(2[(\rho^{-\alpha}-k)^2+k^2] |v|^2\big)(t,x) dx dt\\
&\leq \frac{2}{\mu}(|\alpha+1|^2+2\alpha^2) \int^{t_1}_0\int_{A_k(t)}\big( |v|^2 [(\rho^{-\alpha}-k)^2+k^2] \big) (t,x) dx dt 
\end{aligned}
\label{impcru2}
\end{equation}
Following \cite{La}, we set:
$$|(\rho^{-\alpha})^{(k)}|^2_{Q_{t_1}(k)}=\sup_{0\leq t\leq t_1}\|(\rho^{-\alpha})^{(k)}(t)\|_{L^2}+\|\n (\rho^{-\alpha})^{(k)}\|_{L^2_{t_1}(L^2(\R^N))}.$$
From (\ref{impcru2}) we have by H\"older's inequalities:
\begin{equation}
\min(\frac{1}{2},\frac{\mu}{2})|(\rho^{-\alpha})^{(k)}|^2_{Q_{t_1}(k)}\leq C_{\alpha,\mu}\||v|^2\|_{L^r_{t_1}(L^q(\R^N))}\|(\rho^{-\alpha}-k)^2+k^2\|_{L^{\frac{r}{r-1}}(L^{\frac{q}{q-1}})(Q_{t_1}(k))},
\label{inefond}
\end{equation}
with $Q_{t_1}(k))=\{(t,x); \,t\in [0,t_1], x\in A_k(t)\}$, $C_{\alpha,\mu}=\frac{2}{\mu}(|\alpha+1|^2+2\alpha^2) $ and:
\begin{equation}
\frac{1}{r}+\frac{N}{2q}=1-\kappa_1,
\label{algebrique}
\end{equation}
with $0<\kappa_1<1$. It remains to estimate the right hand side of (\ref{inefond}), let us start with the term $(\rho^{-\alpha}-k)^2$, applying  H\"older's inequality we have:
\begin{equation}
\begin{aligned}
&\|(\rho^{-\alpha}-k)^2\|_{L^{\frac{r}{r-1}}(L^{\frac{q}{q-1}})(Q_{t_1}(k))}=\|(\rho^{-\alpha}-k)\|_{L^{\frac{2r}{r-1}}(L^{\frac{2q}{q-1}})(Q_{t_1}(k))}^2,\\
&\leq \|(\rho^{-\alpha})^{(k)}\|^2_{L^{r_2}_{t_1}(L^{q_2}(\R^N)} |(\int^{t_1}_0  \lambda(A_k(t))^{(\frac{1}{r_1}-\frac{1}{r_2})^{-1}(\frac{1}{q_1}-\frac{1}{q_2})}  dt)^{\frac{1}{r_1}-\frac{1}{r_2}}|^2\\
\end{aligned}
\label{7.9}
\end{equation}
with $q_1=\frac{2q}{q-1}$, $r_1=\frac{2r}{r-1}$, $q_2=q_1(1+\kappa)$,  $r_2=r_1(1+\kappa)$ and $\kappa=\frac{2\kappa_1}{N}$ ( $\lambda$ is the Lebesgue measure). Now we have:
$$\frac{1}{q_1}-\frac{1}{q_2}=\frac{\kappa}{(1+\kappa)q_1},\;\frac{1}{r_1}-\frac{1}{r_2}=\frac{\kappa}{(1+\kappa)r_1}.$$
It implies that using (\ref{7.9}) we have:
\begin{equation}
\begin{aligned}
&\|(\rho^{-\alpha}-k)^2\|_{L^{\frac{r}{r-1}}(L^{\frac{q}{q-1}})(Q_{t_1}(k))}\leq  \|(\rho^{-\alpha})^{(k)}\|^2_{L^{r_2}_{t_ 1}(L^{q_2}(\R^N))}  (\int^{t_1}_0  \lambda(A_k(t))^{\frac{r_1}{q_1}}  dt)^{ \frac{2\kappa}{(1+\kappa)r_1}}\\
\end{aligned}
\label{7.91}
\end{equation}
Using the relation (\ref{algebrique}) and the fact that $\kappa_1=\frac{N\kappa}{2}$ we observe that:
\begin{equation}
\begin{aligned}
\frac{1}{r_2}+\frac{N}{2 q_2}&=\frac{r-1}{2r(1+\kappa)}+\frac{N(q-1)}{4q(1+\kappa)}\\
&=\frac{2+N}{4(1+\kappa)}-\frac{1}{2(1+\kappa)}(\frac{1}{r}+\frac{N}{2q})\\
&=\frac{N}{4(1+\kappa)}+\frac{\kappa_1}{2(1+\kappa)}=\frac{N}{4}.
\end{aligned}
\label{alg2}
\end{equation}
$$q_2\in[2(1+\kappa),\frac{2N(1+\kappa)}{N-2+2\kappa_1}], r_2\in[2(1+\kappa),\frac{2(1+\kappa)}{\kappa_1}],$$
Now using the Gagliardo-Niremberg  inequality we have:
$$ \|(\rho^{-\alpha})^{(k)}(t,\cdot)\|_{L^{q_2}(\R^N)}\leq\beta \|\n (\rho^{-\alpha})^{(k)}(t,\cdot)\|_{L^2(\R^N)}^\alpha\|(\rho^{-\alpha})^{(k)}(t,\cdot)\|_{L^2(\R^N)}^{1-\alpha},$$
with: $\alpha=\N-\frac{N}{q_2}$ and $\alpha=\frac{2}{r_2}$. Using now H\"older's inequality we have:
$$ \|(\rho^{-\alpha})^{(k)}\|_{L^{r_2}_{t_1}(L^{q_2}(\R^N)} \leq \beta\|\n  (\rho^{-\alpha})^{(k)}\|^{1-\frac{2}{r_2}}_{L^2_{t_1}(L^2(\R^N))}\sup_{0\leq t\leq t_1}\|(\rho^{-\alpha})^{(k)}\|_{L^2(\R^N)}^{\frac{2}{r_2}}. $$
Using Young inequality we have:
\begin{equation}
 \|(\rho^{-\alpha})^{(k)}\|_{L^{r_2}_{t_1}(L^{q_2}(\R^N)} \leq\beta |(\rho^{-\alpha})^{(k)}|_{Q_{t_1}}.
 \label{3.42}
 \end{equation}
Plugging this last inequality in (\ref{7.91}) we have:
\begin{equation}
\begin{aligned}
&\|(\rho^{-\alpha}-k)^2\|_{L^{\frac{r}{r-1}}(L^{\frac{q}{q-1}})(Q_{t_1}(k))}\leq \beta^2 |(\rho^{-\alpha})^{(k)}|_{Q_{t_1}}^2\mu(k)^{ \frac{2\kappa}{(1+\kappa)r_1}},
\end{aligned}
\label{7.911}
\end{equation}
with $\mu(k)=  \int^{t_1}_0  \lambda(A_k(t))^{\frac{r_1}{q_1}}  dt$.\\
It remains to bound the term $\|k^2\|_{L^{\frac{r}{r-1}}(L^{\frac{q}{q-1}})(Q_{t_1}(k))}$
\begin{equation}
\begin{aligned}
\|k^2\|_{L^{\frac{r}{r-1}}(L^{\frac{q}{q-1}})(Q_{t_1}(k))}&= k^2\|1\|^2_{L^{r_1}(L^{q_1})(Q_{t_1}(k))},\\
&=\mu(k)^{ \frac{2}{r_1}}.
\end{aligned}
\label{7.11}
\end{equation}
Combining (\ref{7.911}) and (\ref{7.11}) we finally obtain:
\begin{equation}
\begin{aligned}
\min(\frac{1}{2},\frac{\mu}{2}) |(\rho^{-\alpha})^{(k)}|^2_{Q_{t_1}(k)}&\leq C_{\alpha,\mu} \| |v|^2\|_{L^{r}(L^{q})(Q_{t_1}(k))}  \big(\beta^{2} |(\rho^{-\alpha}){(k)}|_{Q_{t_1}(k)}^{2}\mu(k)^{\frac{2\kappa}{(1+\kappa)r_1}}+ k^{2}\mu(k)^{\frac{2}{r_1}}\big),\\
&\leq C_{\alpha,\mu}  \|v\|^2_{L_{t_1}^{2r}(L^{2q}(\R^N))}  
\big(\beta^{2} |(\rho^{-\alpha}){(k)}|_{Q_{t_1}(k)}^{2}\mu(k)^{\frac{2\kappa}{r_1(1+\kappa)}}+ k^{2}\mu(k)^{\frac{2}{r_1}}\big).
\end{aligned}
\end{equation}
with:
$$\frac{1}{r}+\frac{N}{2q}=1-\kappa_1\;\;\mbox{and}\;\;\kappa=\frac{2\kappa_1}{N}.$$
\begin{remarka}
Let us recall that we deal with $t_1\in (0,T^*)$ with $T^*$ the lifespan of our strong solution. It implies in particular that $\frac{1}{\rho}$ belongs in $L^\infty_{t_1}(L^\infty(\R^N))$ for any $t_1$. We are going to prove a postiori that this $L^\infty$ norm depends only of $t_1$ (and not of $T^*$).
\end{remarka}
We deduce that for all $k\geq \sup_{x}\rho_0^{-\alpha}$ and $t_1\in (0,T^*)$:
\begin{equation}
\begin{aligned}
&\min(\frac{1}{2},\frac{\mu}{2})|(\rho^{-\alpha})^{(k)}|^2_{Q_{t_1}(k)}\leq C_{\alpha,\mu}\|\frac{1}{\rho^{\frac{1}{2q}}}\|_{L^\infty_{t_1}(L^\infty)}^2\|\rho^{\frac{1}{2q}}v\|^2_{L^{\infty}_{t_1}(L^{2q}(\R^N))}t_1^{\frac{2}{r}}\\
&\hspace{3cm}\times \big(\beta^{2} |(\rho^{-\alpha}){(k)}|_{Q_{t_1}(k)}^{2}\mu(k)^{\frac{2\kappa}{r_1(1+\kappa)}}+ k^{2}\mu(k)^{\frac{2}{r_1}}\big),\\
&\leq C_{\alpha,\mu} \|\frac{1}{\rho}\|_{L_{t_1}^\infty(L^\infty)}^{\frac{1}{q}}\|\rho^{\frac{1}{2q}}v\|^2_{L^{\infty}_{t_1}(L^{2q}(\R^N))}t_1^{\frac{2}{r}} \big(\beta^{2} |(\rho^{-\alpha}){(k)}|_{Q_{t_1}(k)}^{2}\mu(k)^{\frac{2\kappa}{r_1(1+\kappa)}}+ k^{2}\mu(k)^{\frac{2}{r_1}}\big).
\end{aligned}
\label{4.61im}
\end{equation}
Next we recall that:
$$\mu(k)=\int^{t_1}_0\lambda(A_k(s))^{\frac{r_1}{q_1}}ds.$$
We easily verify that for $\bar{\rho}>1$ and $k\geq 1$:
$$A_k(s)\subset\{x;\,|\sqrt{\rho}-\sqrt{\bar{\rho}}|\geq (\sqrt{\bar{\rho}}-(\frac{1}{k})^{\frac{1}{2\alpha}})\}.$$
And by Tchebytchev inequality we have for $q_3\geq 1$ (that we shall determinate later):
$$\lambda(A_k(s))\leq \frac{\|\sqrt{\rho}(s,\cdot)-\sqrt{\bar{\rho}}\|^{q_3}_{L^{q_3}}}   {(\sqrt{\bar{\rho}}-(\frac{1}{k})^{\frac{1}{2\alpha}})^{q_3}}$$
In particular we deduce that:
\begin{equation}
\mu(k)^{\frac{2\kappa}{r_1(1+\kappa)}}\leq t_1 ^{\frac{2\kappa}{r_1(1+\kappa)}}  \frac{\|\sqrt{\rho}-\sqrt{\bar{\rho}}\|^{\frac{2\kappa q_3}{q_1(1+\kappa)}}_{L_{t_1}^\infty(L^{q_3})}}   {(\sqrt{\bar{\rho}}-(\frac{1}{k})^{\frac{1}{2\alpha}})^{\frac{2\kappa q_3}{q_1(1+\kappa)}  }}
\label{crut1}
\end{equation}
We are now going to choose  $t_1$ such that for any $k\geq \max(1,\sup_x \rho_0^{-\alpha})$:
$$t_1^{\frac{2}{r}+\frac{2\kappa}{r_1(1+\kappa)}}\leq\frac{1}{2}\min(\frac{1}{2},\frac{\mu}{2})\frac{(\sqrt{\bar{\rho}}-(\frac{1}{k})^{\frac{1}{2\alpha}})^{\frac{2\kappa q_3}{q_1(1+\kappa)}  }}{C_{\alpha,\mu} \|\frac{1}{\rho}\|_{L_{t_1}^\infty(L^\infty)}^{\frac{1}{q}}\|\rho^{\frac{1}{2q}}v\|^2_{L^{\infty}_{t_1}(L^{2q}(\R^N))}\beta^2  \|\sqrt{\rho}-\sqrt{\bar{\rho}}\|^{\frac{2\kappa q_3}{q_1(1+\kappa)}}_{L_{t_1}^\infty(L^{q_3})}      } .$$
More precisely since $\frac{2}{r}+\frac{2\kappa}{r_1(1+\kappa)}=\frac{\kappa(r+1)+2}{r(1+\kappa)}$ $t_1$ verifies for $\bar{\rho}>1$: 
\begin{equation}
t_1\leq \biggl( \frac{1}{2}\min(\frac{1}{2},\frac{\mu}{2})\frac{(\sqrt{\bar{\rho}}-1)^{\frac{2\kappa q_3}{q_1(1+\kappa)}  }}{C_{\alpha,\mu} \|\frac{1}{\rho}\|_{L_{t_1}^\infty(L^\infty)}^{\frac{1}{q}}\|\rho^{\frac{1}{2q}}v\|^2_{L^{\infty}_{t_1}(L^{2q}(\R^N))}\beta^2  \|\sqrt{\rho}-\sqrt{\bar{\rho}}\|^{\frac{2\kappa q_3}{q_1(1+\kappa)}}_{L_{t_1}^\infty(L^{q_3})}      } 
 \biggl)^{\frac{r(1+\kappa)}{\kappa(r+1)+2}}.
\label{condit1}
\end{equation}
Combining (\ref{condit1}) and (\ref{4.61im}) it yields:
\begin{equation}
\begin{aligned}
& \frac{1}{2}\min(\frac{1}{2},\frac{\mu}{2})|(\rho^{-\alpha})^{(k)}|^2_{Q_{t_1}(k)}\leq C_{\alpha,\mu} \|\frac{1}{\rho}\|_{L^\infty_{t_1}(L^\infty)}^{\frac{1}{q}}\|\rho^{\frac{1}{2q}}v\|^2_{L^{\infty}_{t_1}(L^{2q}(\R^N))}t_1^{\frac{2}{r}}  k^{2}\mu(k)^{\frac{2}{r_1}}.
\end{aligned}
\end{equation}
We thus get:
\begin{equation}
\begin{aligned}
&|(\rho^{-\alpha})^{(k)}|_{Q_{t_1}(k)}\leq\sqrt{C_{\alpha,\mu}} \|\frac{1}{\rho}\|_{L_{t_1}^\infty(L^\infty)}^{\frac{1}{2q}}\|\rho^{\frac{1}{2q}}v\|_{L^{\infty}_{t_1}(L^{2q}(\R^N))}t_1^{\frac{1}{r}}  k\mu(k)^{\frac{1}{r_1}}.
\end{aligned}
\label{imp1}
\end{equation}
We are going to use the lemma 6.1 from Ladyzenskaya et al p102 \cite{La} that we recall.
\begin{lemma}
 \label{lemLa61}
Assume that $\sup_{x\in\R^N}\frac{1}{\rho_0^{\alpha}(x)}<+\infty$ and that the inequalities
\begin{equation}
|(\rho^{-\alpha})^{(k)}|_{Q_{t_1}(k)}\leq\gamma   k\mu(k)^{\frac{1+\kappa}{r_1(1+\kappa)}}
\label{6.11}
\end{equation}
hold for $k\geq \max(1,\hat{k}_0)$ with $\sup_{x\in\R^N}\frac{1}{\rho_0^{\alpha}(x)}=\hat{k}_0$ with certain positive constants $\gamma$ and $\kappa$. Then we have for $q_3\geq 1$:
\begin{equation}
\begin{aligned}
\sup_{(t,x)\in[0,t_1]\times\R^N}\frac{1}{\rho^{\alpha}(t,x)}&\leq 2\max(1,\sup_{x\in\R^N}\frac{1}{\rho_0^{\alpha}(x)})(1+2^{\frac{2}{\kappa}+\frac{1}{\kappa^2}}(\beta\gamma)^{1+\frac{1}{\kappa}}\\
&\times t_{1}^{\frac{1}{r_1}}(\frac{1}{\sqrt{\bar{\rho}}-1})^{\frac{q_3 }{q_1}}\|\sqrt{\rho}-\sqrt{\bar{\rho}}\|^{\frac{q_3}{q_1}}_{L_{t_1}^{\infty}(L^{q_3})}
\end{aligned}
\label{6.22}
\end{equation}
with $\beta>0$ a constant related with Gagliardo-Niremberg inequality.
\label{lemma4.1}
\end{lemma}
{\bf Proof of Lemma \ref{lemma4.1}:} Let us take the sequence of levels $k_n=M(2-2^{-n})$ with $n\in\mathbb{N}$ and we assume that $M\geq \max(1,\hat{k}_0)>0$. Then we have:
\begin{equation}
(k_{n+1}-k_n)\mu^{r_1(1+\kappa)}(k_{n+1})\leq\| (\rho^{-\alpha})^{(k_n)}\|_{L^{r_2}_{t_1}(L^{q_2}(\R^N))}.
\label{6.3}
\end{equation}
By (\ref{3.42}) and (\ref{6.11}) we have:
\begin{equation}
\| (\rho^{-\alpha})^{(k_n)}\|_{L^{r_2}_{t_1}(L^{q_2}(\R^N))}\leq \beta|(\rho^{-\alpha})^{(k)}|_{Q_{t_1}}\leq \beta\gamma k_{n} \mu(k_n)^{\frac{1+\kappa}{r_1(1+\kappa)}},
\label{6.4}
\end{equation}
and then:
\begin{equation}
\mu^{\frac{1}{r_2}}(k_{n+1})\leq\frac{\beta\gamma k_n}{k_{n+1}-k_{n}}\mu^{\frac{1+\kappa}{r_2}}(k_n)\leq 4\beta\gamma 2^{n}\mu^{\frac{1+\kappa}{r_2}}(k_n). 
\label{suitein}
\end{equation}
Let us recall now the lemma 5.6 p 95 from \cite{La}.
\begin{lemma}
Let $(y_n)_{n\in\mathbb{N}}$ a nonnegative sequence verifying for all $n\geq 0$:
\begin{equation}
y_{n+1}\leq c b^n y_n^{1+\e},
\label{5.14}
\end{equation}
with some positive constants $c$, $\e$ and $b\geq 1$. Then:
\begin{equation}
y_n\leq c^{\frac{(1+\e)^n-1}{\e}}b^{\frac{(1+\e)^n-1}{\e^2}-\frac{n}{\e}}y_0^{(1+\e)^n}.
\label{5.15}
\end{equation}
In particular if:
\begin{equation}
y_0\leq\theta=c^{-\frac{1}{\e}}b^{-\frac{1}{\e^2}}\;\;\mbox{and}\; b>1,
\label{5.16}
\end{equation}
then
\begin{equation}
y_n\leq\theta b^{-\frac{n}{\e}}
\label{5.17}
\end{equation}
and in particular $y_N\rightarrow 0$ for $n\rightarrow +\infty$.
\label{lemma5.6}
\end{lemma}
Applying the lemma \ref{lemma5.6} to the inequality (\ref{suitein}), we deduce that $\mu^{\frac{1}{r_2}}(k_{n})$ goes to zero when $n$ goes to $+\infty$ if  $\mu^{1}{r_2}(k_{0})$ is sufficiently small, namely if:
\begin{equation}
\mu^{\frac{1}{r_2}}(k_{0})=\mu^{\frac{1}{r_2}}(M)\leq(4\beta\gamma )^{-\frac{1}{\kappa}}2^{-\frac{1}{\kappa^2}}.\label{6.6}
\end{equation}
In order to satisfy the previous inequality we set $M=m \hat{k}_0$, $m>1$ and substitute $\hat{k}$  for $k_n$ and $M$ for $k_{n+1}$ in (\ref{6.3}) this gives by Tchebytchev inequality:
$$
\begin{aligned}
&\mu^{\frac{1}{r_2}}(M)\leq\frac{\beta\gamma}{m-1}\mu^{\frac{1}{r_1}}(\max(1,\hat{k}_0))\leq 
\frac{\beta\gamma}{m-1} t_1 ^{\frac{1}{r_1}}  \frac{\|\sqrt{\rho}-\sqrt{\bar{\rho}}\|^{\frac{q_3}{q_1}}_{L_{t_1}^\infty(L^{q_3})}}   {(\sqrt{\bar{\rho}}-(\frac{1}{\max(1,\hat{k}_0))})^{\frac{1}{2\alpha}})^{\frac{ q_3}{q_1}  }}.
\end{aligned}
$$
It suffice to choose:
$$m=1+\beta\gamma  t_1 ^{\frac{1}{r_1}}\frac{\|\sqrt{\rho}-\sqrt{\bar{\rho}}\|^{\frac{q_3}{q_1}}_{L_{t_1}^\infty(L^{q_3})}}   {(\sqrt{\bar{\rho}}-(\frac{1}{\max(1,\hat{k}_0))})^{\frac{1}{2\alpha}})^{\frac{ q_3}{q_1}  }}
(4\beta\gamma)^{\frac{1}{\kappa}}2^{\frac{1}{\kappa^2}}.$$
Then the condition (\ref{6.6}) is verified and then $\mu(2M)$ is equal to zero which implies that:
$$\|\rho^{-\alpha}\|_{L^\infty([0,t_1]\times\R^N)}\leq 2M=2m\max(1,\hat{k}_0).$$
$\blacksquare$\\
\\
We deduce that choosing $t_1$ as in (\ref{condit1}) we have via the lemma \ref{lemLa61} and the fact that $\gamma=C_{\alpha,\mu} \|\frac{1}{\rho}\|_{L_{t_1}^\infty(L^\infty)}^{\frac{1}{2q}}\|\rho^{\frac{1}{2q}}v\|_{L^{\infty}_{t_1}(L^{2q}(\R^N))}t_1^{\frac{1}{r}} $:
\begin{equation}
\begin{aligned}
\sup_{(t,x)\in[0,t_1]\times\R^N}\frac{1}{\rho^{\alpha}(t,x)}&\leq 2\sup_{x\in\R^N}\frac{1}{\rho_0^{\alpha}(x)}\biggl(1+2^{\frac{2}{\kappa}+\frac{1}{\kappa^2}}(\beta  C_{\alpha,\mu} \|\frac{1}{\rho}\|_{L^\infty(L^\infty)}^{\frac{1}{2q}}\|\rho^{\frac{1}{2q}}v\|_{L^{\infty}_{t_1}(L^{2q}(\R^N))}t_1^{\frac{1}{r}}   )^{1+\frac{1}{\kappa}}t_1^{\frac{1}{r_1}}\\
&\times (\frac{1}{\sqrt{\bar{\rho}}-1})^{\frac{q_3}{q_1}}\|\sqrt{\rho}-\sqrt{\bar{\rho}}\|^{\frac{q_3}{q_1}}_{L_{t_1}^{\infty}(L^{q_3})}\biggl).
\end{aligned}
\label{cru3}
\end{equation}
In order to simplify the notations we are going to deal with $\alpha=1$, 
$\kappa_1=\frac{1}{2}$ such that $\kappa=\frac{1}{N}$ and:
\begin{equation}
\frac{1}{r}+\frac{N}{2 q}=\frac{1}{2}.
\label{rela}
\end{equation}
We have then since $\sup_{(t,x)\in[0,t_1]\times\R^N}\frac{1}{\rho}(t,x)=\|\frac{1}{\rho}\|_{L^{\infty}_{t_1}(L^{\infty})}$ and setting abusively $C_{\alpha,\mu}\beta=\beta$ (in order to simplify the notation):
\begin{equation}
\begin{aligned}
 \|\frac{1}{\rho}\|_{L^\infty_{t_1}(L^\infty)}  &\leq 2  \max(1, \|\frac{1}{\rho_0}\|_{L^\infty})  \biggl(1+  \|\frac{1}{\rho}\|_{L^\infty_{t_1}(L^\infty)}^{\frac{1}{2q} (N+1)}2^{2N+N^2}(\beta\|\rho^{\frac{1}{2q}}v\|_{L^{\infty}_{t_1}(L^{2q}(\R^N))}t_1^{\frac{1}{r}}   )^{N+1}t_1^{\frac{1}{r_1}}\\
&\times (\frac{1}{\sqrt{\bar{\rho}}-1})^{\frac{q_3}{q_1}}\|\sqrt{\rho}-\sqrt{\bar{\rho}}\|^{\frac{q_3}{q_1}}_{L_{t_1}^{\infty }(L^{q_3})}\biggl).
\end{aligned}
\label{cru4}
\end{equation}
Choosing $q$ large enough and applying the Young inequality with $p=\frac{2q}{N+1}$ and $p'=\frac{2q}{2q-N-1}$ we have:
\begin{equation}
\begin{aligned}
 \|\frac{1}{\rho}\|_{L^\infty_{t_1}(L^\infty)}  &\leq \frac{N+1}{2q} \|\frac{1}{\rho}\|_{L_{t_1}^\infty(L^\infty)}+ \frac{2q-N-1}{2q} \biggl(2\max(1, \|\frac{1}{\rho_0}\|_{L^\infty})2^{2N+N^2}(\beta\|\rho^{\frac{1}{2q}}v\|_{L^{\infty}_{t_1}(L^{2q}(\R^N))}t_1^{\frac{1}{r}}   )^{N+1} \\
&\times t_1^{\frac{1}{r_1}} (\frac{1}{\sqrt{\bar{\rho}}-1})^{\frac{q_3}{q_1}}\|\sqrt{\rho}-\sqrt{\bar{\rho}}\|^{\frac{q_3}{q_1}}_{L_{t_1}^{\infty }(L^{q_3})}\biggl)^{\frac{2q}{2q-N-1}}+ 2  \max(1, \|\frac{1}{\rho_0}\|_{L^\infty}) .
\end{aligned}
\label{cru5}
\end{equation}
It gives:
\begin{equation}
\begin{aligned}
 \|\frac{1}{\rho}\|_{L^\infty_{t_1}(L^\infty)}  &\leq  \max(1,\|\frac{1}{\rho_0}\|_{L^\infty} )\biggl(\frac{4q}{2q-N-1}+\max(1, \|\frac{1}{\rho_0}\|_{L^\infty})^{\frac{N+1}{2q-N-1}}\big(2 2^{2N+N^2}\\
 &\times(\beta\|\rho^{\frac{1}{2q}}v\|_{L^{\infty}_{t_1}(L^{2q}(\R^N))}t_1^{\frac{1}{r}}   )^{N+1}t_1^{\frac{1}{r_1}}  (\frac{1}{\sqrt{\bar{\rho}}-1})^{\frac{q_3}{q_1}}\|\sqrt{\rho}-\sqrt{\bar{\rho}}\|^{\frac{q_3}{q_1}}_{L_{t_1}^{\infty}(L^{q_3})}\big)^{\frac{2q}{2q-N-1}}\biggl).
\end{aligned}
\label{cru6}
\end{equation}
It implies in particular: 
\begin{equation}
\begin{aligned}
 \|\frac{1}{\rho}\|_{L^\infty_{t_1}(L^\infty)}  &\leq \max(1, \|\frac{1}{\rho_0}\|_{L^\infty})^{1+\frac{N+1}{2q-N-1}} \biggl(2+\big(2 2^{2N+N^2}(\beta\|\rho^{\frac{1}{2q}}v\|_{L^{\infty}_{t_1}(L^{2q}(\R^N))}t_1^{\frac{1}{r}}   )^{N+1}t_1^{\frac{1}{r_1}}\\
&\times (\frac{1}{\sqrt{\bar{\rho}}-1})^{\frac{q_3}{q_1}}\|\sqrt{\rho}-\sqrt{\bar{\rho}}\|^{\frac{q_3 }{q_1}}_{L_{t_1}^{\infty}(L^{q_3})}\big)^{\frac{2q}{2q-N-1}}\biggl).
\end{aligned}
\label{cru7}
\end{equation}
Using the lemma \ref{gainin} and the proposition \ref{22energy}  with $q_3=2$, we deduce that for $t_1\in (0,T^*)$ and $t_1$ verifying the assumption (\ref{condit1}) we have for any large $q$:
\begin{equation}
\begin{aligned}
 \|\frac{1}{\rho}\|_{L^\infty_{t_1}(L^\infty)}  &\leq \max(1, \|\frac{1}{\rho_0}\|_{L^\infty})^{1+\frac{N+1}{2q-N-1}} \biggl(2+\big(2 2^{2N+N^2}(\beta C(T^{*})\times(T^{*})^{\frac{1}{r}}   )^{N+1}(T^*)^{\frac{1}{r_1}}\\
&\times ((\frac{1}{\sqrt{\bar{\rho}}-1})^{\frac{q-1}{q}}E_0^{\frac{q-1 }{q}}\big)^{\frac{2q}{2q-N-1}}\biggl)\\
&\leq  \max(1, \|\frac{1}{\rho_0}\|_{L^\infty})^{1+\frac{N+1}{2q-N-1}} \biggl(2+\big(2 2^{2N+N^2}(\beta C(T^{*})\times(T^{*})^{\frac{1}{r}}   )^{N+1}(T^*)^{\frac{1}{r_1}}\big)^{\frac{2q}{2q-N-1}}   \\
&   \times\big((\frac{1}{\sqrt{\bar{\rho}}-1})E_0\big)^{\frac{2(q-1)}{2q-N-1}}\biggl),
\end{aligned}
\label{cru7b}
\end{equation}
with $E_0$ depending only on the initial data. It implies that for any $q$ it exists $C$ large enough and $\gamma>0$ depending only on the norm of the initial data $\|\rho_0^{\frac{1}{q}}v\|_{L^q}$ and from $E_0$ such that:
\begin{equation}
\begin{aligned}
 \|\frac{1}{\rho}\|_{L^\infty_{t_1}(L^\infty)}  &\leq \max(1, \|\frac{1}{\rho_0}\|_{L^\infty})^{1+\frac{N+1}{2q-N-1}} (2+C(1+(T^*)^{\gamma})).
\end{aligned}
\label{cru7d}
\end{equation}
\begin{remarka}
Since with the choice on the initial data in the theorem \ref{theo2}, $\|\rho_0^{\frac{1}{q}}v\|_{L^q}$ is uniformly bounded in $q$, it implies that $C$ and $\gamma$ do not depend on $q$.
 \end{remarka}
Combining the condition (\ref{condit1}) and the previous inequality (\ref{cru7d}) $t_1$ must verify the following condition with $C_{\alpha,\mu}=C_\mu$ by notation:
\begin{equation}
t_1\leq \biggl( \frac{1}{2}\min(\frac{1}{2},\frac{\mu}{2})\frac{(\sqrt{\bar{\rho}}-1)^{\frac{2\kappa (q-1)}{q(1+\kappa)}  }}{ C_\mu \max(1, \|\frac{1}{\rho_0}\|_{L^\infty})^{\frac{2}{2q-N-1}} (2+C(1+(T^*)^{\gamma}))^{\frac{1}{q}}\|\rho^{\frac{1}{2q}}v\|^2_{L^{\infty}_{t_1}(L^{2q})}\beta^2  \|\sqrt{\rho}-\sqrt{\bar{\rho}}\|^{\frac{2\kappa (q-1)}{q(1+\kappa)}}_{L^\infty_{t_1}(L^{2})}      } \biggl)^{\frac{r(1+\kappa)}{\kappa(r+1)+2}}.
\label{t1cru}
\end{equation}
We deduce using lemma \ref{gainin} and proposition \ref{22energy} that for $q\geq 1$:
\begin{equation}
\begin{aligned}
&t_1\leq \frac{M}{M'(1+T^{\beta})}\frac{1}{ \max(1, \|\frac{1}{\rho_0}\|_{L^\infty})^{\frac{2}{2q-N-1}\alpha_1 }},\\
& \|\frac{1}{\rho}\|_{L^\infty_{t_1}(L^\infty)}  \leq \max(1, \|\frac{1}{\rho_0}\|_{L^\infty})^{1+\frac{N+1}{2q-N-1}} (2+C(1+(T^*)^{\gamma})).
\end{aligned}
\label{estimt1cru}
\end{equation}
with  $\alpha_1=\frac{r(1+\kappa)}{\kappa(r+1)+2}$. Here $C>0$, $M,M'>0$, $\beta>0$ depends only on the initial data $\|\rho_0^{\frac{1}{q}}v_0\|_{L^q}$, $E_0$ and is independent of $q$ and of $t_1$. In an other way we have:
\begin{equation}
\begin{aligned}
&t_1\leq \frac{C_1}{ \max(1, \|\frac{1}{\rho_0}\|_{L^\infty})^{\frac{2}{2q-N-1}\alpha_1 }},\\
& \|\frac{1}{\rho}\|_{L^\infty_{t_1}(L^\infty)}  \leq C_2\max(1, \|\frac{1}{\rho_0}\|_{L^\infty})^{1+\frac{N+1}{2q-N-1}}.
\end{aligned}
\label{estimt1cru2}
\end{equation}
with $C_1=\frac{M}{M'(1+T^{\beta})}$ and $C_2= (2+C(1+(T^*)^{\gamma}))$.\\
For the moment we have proved estimates on $\frac{1}{\rho}$ in $L^\infty$ norm when $t_1$ verifies the condition (\ref{estimt1cru2}), we would like to get the same type of estimate for any $T\in (0,T^*)$. To do this we are going to repeat the argument used for the case $t_1$, it implies that by
recurrence  it exists a sequel $(t_n)_{n\in\mathbb{N}}$ verifying for $n\geq 2$:
\begin{equation}
\begin{aligned}
&t_n\leq \frac{C_1}{ C^{(\sum_{i=0}^{n-2}\beta^{i}_1 )\frac{2}{2q-N-1}\alpha_1}_2\max(1, \|\frac{1}{\rho_0}\|_{L^\infty})^{\beta_1^{n-1}\frac{2}{2q-N-1}\alpha_1 }},\\
& \|\frac{1}{\rho}\|_{L^\infty_{\sum_{i=1}^{n}t_i}(L^\infty)}  \leq C_2^{\sum_{i=0}^{n-1}\beta_1^{i}} \max(1, \|\frac{1}{\rho_0}\|_{L^\infty})^{\beta_1^n}.
\end{aligned}
\label{estimt1cru3}
\end{equation}
with $\beta_1=1+\frac{N+1}{2q-N-1}$.\\
Our goal is now to verify that  $\sum_{i=1}^{+\infty}t_i$ can cover any $T\in (0,T^*)$ or in other terms that  $\sum_{i=1}^{+\infty}t_i> T$  for any $T\in (0,T^*)$ when $q$ is large enough. It will prove that it exists $n$ large enough such that $\sum_{i=1}^{n-1}t_i\leq T<\sum_{i=1}^{n}t_i$ and this for any $T\in (0,T^*)$ (with $n$ depending on $T$). We have then:
\begin{equation}
\sum_{i=1}^{n}t_i=\frac{C_1}{ \max(1, \|\frac{1}{\rho_0}\|_{L^\infty})^{\frac{2}{2q-N-1}\alpha_1 }}+\sum_{j=2}^{n} \frac{C_1}{ C^{(\sum_{i=0}^{j-2}\beta^{i}_1 )\frac{2}{2q-N-1}\alpha_1}_2\max(1, \|\frac{1}{\rho_0}\|_{L^\infty})^{\beta_1^{j-1}\frac{2}{2q-N-1}\alpha_1 }}
\label{calcultn}
\end{equation}
Next we have:
$$\sum_{i=0}^{j-2}\beta^{i}_1=\frac{1-\beta_1^{j-1}}{1-\beta_1}= \frac{2q-N-1}{N+1}\big((1+\frac{N+1}{2q-N-1})^{j-1}-1\big)     .$$
It implies that:
\begin{equation}
(\sum_{i=0}^{j-2}\beta^{i}_1 )\frac{2}{2q-N-1}\alpha_1=\frac{2\alpha_1}{N+1}\big((1+\frac{N+1}{2q-N-1})^{j-1}-1\big) 
\end{equation}
Now we are going to consider the $\lim_{q\rightarrow +\infty} \sum_{i=1}^{n}t_i$ in order to prove that this sum converges to $+\infty$. We have then:
\begin{equation}
\begin{aligned}
&(\sum_{i=0}^{j-2}\beta^{i}_1 )\frac{2}{2q-N-1}\alpha_1\sim_{q\rightarrow+\infty} (j-1)\frac{2\alpha_1}{2q-N-1}\\
&\beta_1^{j-1}\frac{2\alpha_1}{2q-N-1}\sim_{q\rightarrow+\infty} \frac{2\alpha_1}{2q-N-1}
\end{aligned}
\end{equation}
Here we recall via the relation (\ref{rela}) that $r$ goes to $2$ when $q$ goes to $+\infty$. It implies that $\alpha_1$ converges to $\frac{2(1+\frac{1}{N})}{\frac{3}{N}+2}$.
We deduce via the previous estimates and the equality (\ref{calcultn}) that:
\begin{equation}
\lim_{q\rightarrow+\infty}\sum_{i=1}^{n}t_i=C_1 n.
\end{equation}
And for $n$ large enough depending on $T$, $q$ large enough we have the existence of $n$ such that $\sum_{i=1}^{n-1}t_i\leq T<\sum_{i=1}^{n}t_i$. It concludes the proof of the proposition \ref{avide}
since via the estimate (\ref{estimt1cru3}) we observe that $\|\frac{1}{\rho}\|_{L^\infty_T(L^\infty(\R^N))}$ depends only on $n$, $q$ and  $\|\frac{1}{\rho_0}\|_{L^\infty}$ but $n$ depends on $T$. This is exactly what we want. $\blacksquare$
\subsection{Control on $\|(q,v)(T,\cdot)\|_{B^{\NN+\e'}_{p,1}\times B^{\NN-1+\e'}_{p,1}}$ for any $T<T^{*}$ and $\e'>0$}
\subsubsection{Estimate on $\|v(T,\cdot)\|_{ B^{\NN-1+\e'}_{p,1}}$ for any $T<T^{*}$ and $\e'>0$}
Let us recall that we have proved in the lemma \ref{gainin} and proposition \ref{avide} that for any $p\geq 2$ it exists increasing functions $C$ and $C_1$ with $C$ and $C_1$ depending only on the initial data $(\rho_0,v_0)$ such that:
\begin{equation}
\begin{aligned}
&\|\rho^{\frac{1}{p}}v\|_{L_T^{\infty}(L^p)}\leq C(T)\\
&\|\frac{1}{\rho^{\frac{1}{p}} }\|_{L_T^{\infty}(L^{\infty})}\leq C_1^{\frac{1}{p}}(T).
\end{aligned}
\end{equation}
It yields:
\begin{equation}
\|v\|_{L^{\infty}_T(L^p)}\leq C(T)C_1^{\frac{1}{p}}(T).
\end{equation}
By Besov embedding we observe that for $p\geq 2$ it exists $C>0$ such that:
\begin{equation}
\begin{aligned}
&\|v\|_{L^{\infty}_T(B^{0}_{p,\infty})}\leq CC(T)C_1^{\frac{1}{p}}(T),\\
&\|v\|_{L^{\infty}_T(B^{-N(\frac{1}{2}-\frac{1}{p})}_{p,\infty})}\leq C C(T)C_1^{\frac{1}{2}}(T).
\end{aligned}
\label{inclusion}
\end{equation}
By interpolation in the Besov space, we deduce that when $-\N+\NN<\NN-1+\e'<0$ (it corresponds to $p>\frac{N}{1-\e'}$) we have for $M$ an increasing function depending only on the initial data $(\rho_0,v_0)$:
\begin{equation}
\begin{aligned}
&\|v\|_{\widetilde{L}^{\infty}_T(B^{\NN-1+\e'}_{p,1})} \leq M(T).
\end{aligned}
\label{inclusion}
\end{equation}
Since we know that $v$ belongs in $\widetilde{C}_{T^*}(B^{\NN-1+\e'}_{p,1})$ (\ref{inclusion}) implies that:
\begin{equation}
\|v(T,\cdot)\|_{B^{\NN-1+\e'}_{p,1}}\leq M(T).
\label{4.98}
\end{equation}
$\blacksquare$
\subsubsection{Estimate on $\|q(T,\cdot)\|_{ B^{\NN+\e'}_{p,1}}$ for any $T<T^{*}$ and $\e'>0$}
It suffices to use the first equation in (\ref{systK}) and the proposition \ref{chaleur} on the heat equation which ensures that for any $T\in(0,T^*)$ we have for $C>0$ and $q_1=\rho-\bar{\rho}$:
\begin{equation}
\begin{aligned}
\|q_1\|_{\widetilde{L}^\infty_T(B^{\NN+\e'}_{p,1})}&\leq C (\|q^1_0\|_{B^{\NN+\e'}_{p,1}}+\|{\rm div}(\rho v)\|_{\widetilde{L}^\infty_T(B^{\NN-2+\e'}_{p,1})})\\
&\leq C (\|q^1_0\|_{B^{\NN+\e'}_{p,1}}+\|\rho v\|_{\widetilde{L}^\infty_T(B^{\NN-1+\e'}_{p,1})}).
\end{aligned}
\label{impoq2}
\end{equation}
Proceeding as in the previous section we have for any $p\geq 2$:
\begin{equation}
\begin{aligned}
&\|\rho v\|_{\widetilde{L}^\infty_T(B^{0}_{p,\infty})}=\|\rho v\|_{L^\infty_T(B^{0}_{p,\infty})}\leq \|\rho\|^{1-\frac{1}{p}}_{L^\infty_T(L^\infty)} \|\rho^{\frac{1}{p}} v\|_{L^\infty_T(L^p)},\\
&\|\rho v\|_{\widetilde{L}^\infty_T(B^{-N(\frac{1}{2}-\frac{1}{p})}_{p,\infty})}=\|\rho v\|_{L^\infty_T(B^{-N(\frac{1}{2}-\frac{1}{p})}_{p,\infty})}\leq \|\rho\|^{\frac{1}{2}}_{L^\infty_T(L^\infty)} \|\rho^{\frac{1}{2}} v\|_{L^\infty_T(L^2)}.
\end{aligned}
\label{impov2}
\end{equation}
By interpolation we deduce that for $-\N+\NN<\NN-1+\e'<0$ we have with $\NN-1+\e'=\theta(-\N+\NN)$ and $(1-\theta)(1-\frac{1}{p})+\frac{\theta}{2}=1-\frac{1-\e'}{N}$:
\begin{equation}
\begin{aligned}
\|\rho v\|_{\widetilde{L}^\infty_T(B^{\NN-1+\e'}_{p,1})}&\leq \|\rho\|^{\frac{\theta}{2}}_{L^\infty_T(L^\infty)} \|\rho^{\frac{1}{2}} v\|^\theta_{L^\infty_T(L^2)} \|\rho\|^{(1-\theta)(1-\frac{1}{p})}_{L^\infty_T(L^\infty)} \|\rho^{\frac{1}{p}} v\|^{1-\theta}_{L^\infty_T(L^p)}\\
&\leq  \|\rho\|^{1-\frac{1-\e'}{N}}_{L^\infty_T(L^\infty)} \|\rho^{\frac{1}{2}} v\|^\theta_{L^\infty_T(L^2)} \|\rho^{\frac{1}{p}} v\|^{1-\theta}_{L^\infty_T(L^p)}.
\end{aligned}
\label{estimv1}
\end{equation}
Using the lemme \ref{gainin} and (\ref{estimv1}) it exists a increasing function $M_1$ such that:
\begin{equation}
\begin{aligned}
\|\rho v\|_{\widetilde{L}^\infty_T(B^{\NN-1+\e'}_{p,1})} \leq \|\rho\|^{1-\frac{1-\e'}{N}}_{L^\infty_T(L^\infty)}M_1(T).
\end{aligned}
\label{estimv1}
\end{equation}
Plugging the previous estimate in (\ref{impov2}) we have:
\begin{equation}
\begin{aligned}
\|q_1\|_{\widetilde{L}^\infty_T(B^{\NN+\e'}_{p,1})}
&\leq C (\|q^1_0\|_{B^{\NN+\e'}_{p,1}}+(\|q_1\|_{L^\infty_T(L^\infty)} +\bar{\rho})^{1-\frac{1-\e'}{N}}M_1(T)).
\end{aligned}
\label{impoq2}
\end{equation}
Now by Besov embedding and interpolation, we recall that it exists $C>0$ such that:
\begin{equation}
\begin{aligned}
\|q_1\|_{L^\infty_T(L^\infty)}&\leq C\|q_1\|_{\widetilde{L}^\infty_T(B^{\NN}_{p,1})},\\
&\leq   C'\|q_1\|^\theta_{\widetilde{L}^\infty_T(B^{-N(\frac{1}{2}-\frac{1}{p})}_{p,1})} \|q_1\|^{1-\theta}_{\widetilde{L}^\infty_T(B^{\NN+\e'}_{p,1})},
\end{aligned}
\label{controleL}
\end{equation}
with $\NN=-\theta N(\frac{1}{2}-\frac{1}{p})+(1-\theta)(\NN+\e')$. Next since via the lemma \ref{22energy} we know that $\sqrt{\rho}-\sqrt{\bar{\rho}}$ is in $L^\infty_T(H^1(\R^N))$. We deduce in particular that $q_1=(\rho-\bar{\rho})$ is in $L^\infty_T(L^2(\R^N))$ when $N\leq 4$. Indeed we have:
$$\rho-\bar{\rho}=(\sqrt{\rho}-\sqrt{\bar{\rho}})^2+2\sqrt{\bar{\rho}}(\sqrt{\rho}-\sqrt{\bar{\rho}}),$$
and we conclude using the fact that $(\sqrt{\rho}-\sqrt{\bar{\rho}})^2$ is in $L^\infty_T(L^2(\R^N))$ by Sobolev embedding when $N\leq 4$. When $N\geq 5$ it suffices to observe that $q_1$ is in $L¬\infty_T(L^1(\R^N)+L^2(\R^N))$ and using again interpolation. By Young inequality and the lemma \ref{22energy} we deduce that it exists $M'$ an increasing function such that:
\begin{equation}
\begin{aligned}
\|q_1\|_{L^\infty_T(L^\infty)}\leq &M'(T) +\|q_1\|_{\widetilde{L}^\infty_T(B^{\NN+\e'}_{p,1})},
\end{aligned}
\label{controleL1}
\end{equation}
By Young inequality and Besov embedding we have for $C>0$ and any $\e>0$:
\begin{equation}
\begin{aligned}
&\|q_1\|_{\widetilde{L}^\infty_T(B^{\NN+\e'}_{p,1})}
\leq C (\|q^1_0\|_{B^{\NN+\e'}_{p,1}}+(\|q_1\|_{\widetilde{L}^\infty_T(B^{\NN+\e'}_{p,1})} +M'(T))^{1-\frac{1-\e'}{N}}M_1(T))\\
&\leq C (\|q^1_0\|_{B^{\NN+\e'}_{p,1}}+\frac{\e(N-1+\e')}{N}(\|q_1\|_{\widetilde{L}^\infty_T(B^{\NN+\e'}_{p,1})} +M'(T))+\frac{1-\e'}{\e N} M_1(T)^{\frac{N}{1_\e'}}).
\end{aligned}
\label{impoq3}
\end{equation}
Choosing $\e$ sufficiently small we deduce that it exists $C>0$ such that:
\begin{equation}
\begin{aligned}
\|q_1\|_{\widetilde{L}^\infty_T(B^{\NN+\e'}_{p,1})}
&\leq C (1+\|q^1_0\|_{B^{\NN+\e'}_{p,1}}+M'(T)+ M_1(T)^{\frac{N}{1-\e'}})
\end{aligned}
\label{impoq4}
\end{equation}
It implies in particular that $\rho$ is in $L^\infty_T(L^\infty)$. Now since $q=\ln(\frac{\rho}{\bar{\rho}})=\ln (\frac{q_1+\bar{\rho}}{\bar{\rho}})$ we deduce by proposition \ref{composition}, (\ref{impoq4}) and 
the fact that $\rho$ and $\frac{1}{\rho}$ are in $L^\infty_T(L^\infty)$ (see the proposition \ref{avide}) that it exists $M_2$, $M_3$ increasing function and $C>0$ such that:
\begin{equation}
\begin{aligned}
\|q\|_{\widetilde{L}^\infty_T(B^{\NN+\e'}_{p,1})}
&\leq C (1+\|q^1_0\|_{B^{\NN+\e'}_{p,1}}+M_3(T)+ M_2(T)^{\frac{N}{1-\e'}}).
\end{aligned}
\label{impoq4}
\end{equation}
In particular since $q$ belongs in $\widetilde{C}_{T^*}(B^{\NN+\e'}_{p,1})$ we have for any $T\in (0,T^*)$:
\begin{equation}
\|q(T,\cdot)\|_{B^{\NN+\e'}_{p,1}}
\leq C (1+\|q^1_0\|_{B^{\NN+\e'}_{p,1}}+M_3(T)+ M_2(T)^{\frac{N}{1-\e'}}).
\label{estimq7}
\end{equation}
\section{Proof of the theorem \ref{theo2}}
\label{section5}
Since we have assumed that $T^*<+\infty$, we are interested in proving that this is absurd. It suffices to 
extend the strong solution $(q,v)$ of the system (\ref{systK1}) beyond $T^*$.  Let us summarize which estimates we have obtained on $q$, for all $T\in(0,T^*)$ and $p>\frac{N}{1-\e'}$ with $N\geq 2$ we have via (\ref{4.98}) and (\ref{estimq7}) :
\begin{equation}
\begin{aligned}
&\|v(T,\cdot)\|_{B^{\N-1+\e'}_{p,1}}\leq M(T),\\
&\|q(T,\cdot)\|_{B^{\N+\e'}_{p,1}}\leq M_3(T),
\end{aligned}
\end{equation}
with $M$ and $M_3$ increasing function depending only on the initial data $(q_0,v_0)$. Let us start with studying the case $N\geq 3$. By proceeding as in the previous section we who that:
\begin{equation}
\begin{aligned}
&\|v(T,\cdot)\|_{B^{\N-1}_{p,1}}\leq M_4(T),\\
&\|q(T,\cdot)\|_{B^{\N}_{p,1}}\leq M_5(T),
\end{aligned}
\end{equation}
with $M_4$ and $M_5$ increasing function depending only on the initial data $(q_0,v_0)$ (indeed we need to verify the condition $-\N+\N<\N-1<0$).\\
 Let us consider now the system (\ref{systK1}) with initial data $(q(T^*-\e_1,\cdot),v(T^*-\e_1,\cdot))$ such that:
\begin{equation}
\e_1< D(M(T),M_3(T)).
\label{argfin}
\end{equation}
Using the theorem \ref{theo1} we know that it exists a strong solution $(q_1,v_1)$ to the system (\ref{systK1}) with initial data $(q(T^*-\e_1,\cdot),v(T^*-\e_1,\cdot))$ such that the time of existence $T^*_1$ verifies via the estimate (\ref{3.34}):
$$T^*_1\geq D(M(T),M_3(T))>\e_1.$$
It implies that $T^*-e_1+T^*_1>T^*$, in addition using the uniqueness part of the theorem \ref{theo1} we observe that:
$$(q_{1}(t,\cdot),v_{1}(t,\cdot))=(q(T^*-\e_1+t,\cdot),v(T^*-\e_1+t,\cdot))\;\;\mbox{on}\;(0,\e_1).$$
It implies that we can extend the solution $(q,v)$ beyond $T^*$ which implies that $T^*<+\infty$ is absurd. We have proved then that $T^*=+\infty$ when $N\geq 3$.\\
Let us deal now with the case $N=2$. The difficulty here corresponds to estimate $v$ in $\widetilde{L}_T^\infty(B^{\NN-1}_{p,1})$ using interpolation argument. The problem comes from the low frequencies since $v$ belongs in $L^\infty_T(L^2)$ using the energy estimate and the fact that $\frac{1}{\rho}$ is in $L^\infty_T(L^\infty)$. It implies that $v$ is only in $L^\infty_T(B^{\N-1}_{2,2})$ which is embedded in $L^\infty_T(B^{\NN-1}_{p,2})$ when $p\geq 2$. We have no hope to estimate $\|v(T,\cdot)\|_{B^{\NN-1}_{p,1}}$, we need to proceed in a different way.\\
It suffices to repeat all the procedure with initial data $(q_0,v_0)$ which belong in $B^{\NN}_{p,2}\times B^{\NN-1}_{p,2}$. We can prove the existence of strong solution in finite time (we refer to \cite{Hprepa}, let us mention that to prove this result the choice of the physical coefficient is crucial, in particular the fact that $P(\rho)=a\rho$. It allows to avoid any control on the $L^\infty$ norm of the density what is important in general in order to apply composition theorem). Newt in a similar way we bounded by below the time of existence. And using the same previous arguments we obtain that for $\e'>0$:
\begin{equation}
\begin{aligned}
&\|v(T,\cdot)\|_{B^{\N-1+\e'}_{p,2}}\leq M'(T),\\
&\|q(T,\cdot)\|_{B^{\N+\e'}_{p,2}}\leq M'_3(T),
\end{aligned}
\end{equation}
with $M'$ and $M'_3$ increasing function depending only on the initial data $(q_0,v_0)$ and:
\begin{equation}
\begin{aligned}
&\|v(T,\cdot)\|_{B^{\N-1}_{p,2}}\leq M'_4(T),\\
&\|q(T,\cdot)\|_{B^{\N}_{p,2}}\leq M'_5(T),
\end{aligned}
\end{equation}
with $M'_4$ and $M'_5$ increasing function depending only on the initial data $(q_0,v_0)$.
The conclusion is now the same that for the case $N\geq 3$ and we have $T^*=+\infty$ when $N=2$. It concludes the proof of the theorem \ref{theo2}.
$\blacksquare$
\section{Proof of the theorem \ref{theo2} for general pressure}
\label{section6}
Let us start by giving some basic information on the energy setimate when $P(\rho)=a\rho^\gamma$ with $\gamma>1$. We can check then that $\rho-1$ belongs to the Orlicz space  $L^\infty_T(L^\gamma_2(\R^N))$ (se \cite{fL2}) via the energy estimate.\\
The proof now follows exactly the same line than  the particular case $P(\rho)=a\rho$ with one exception, indeed we have to proceed in an other way if we want to obtain uniform estimate in $p$  on $\rho^{\frac{1}{p}}v$ in $L^\infty_{loc}(L^p(\R^N))$  (for any $p\geq 2$). Following (\ref{in1A1}), we have for any  $0<t\leq T^*$ and $p\geq 4$:
\begin{equation}
\begin{aligned}
&\frac{1}{p}\int_{\R^{N}}(\rho|v|^{p})(t,x)dx+\int^{t}_{0}\int_{\R^{N}} \rho|v|^{p-2}|\n v|^{2}(s,x)dsdx\\
&+\frac{(p-2)}{4}\int^{t}_{0} \int_{\R^{N}}\rho|\n(|v|^{2})|^{2}|v|^{p-4}(s,x)dsdx\leq \frac{1}{p}\int_{\R^{N}}(\rho_{0}|v_{0}|^{p})(x)dx\\
&\hspace{7cm}+|\int^{t}_{0} \int_{\R^{N}}|v|^{p-2}v\cdot\n a\rho^\gamma(s,x)dsdx|.
\end{aligned}
\label{in1A1bn}
\end{equation}
\subsection{Case $N=2$}
Let us start with the case $N=2$. We have essentially to deal with the term:
\begin{equation}
\begin{aligned}
&|\int^{t}_{0} \int_{\R^{N}}|v|^{p-2}v\cdot\n a\rho^\gamma(s,x)dsdx|.
\end{aligned}
\label{pression}
\end{equation}
We have by H\"older and Young inequality:
\begin{equation}
\begin{aligned}
&|\int^{t}_{0} \int_{\R^{N}}|v|^{p-2}v\cdot\n a\rho^\gamma(s,x)dsdx|\leq C |\int^{t}_{0} \int_{\R^{N}}\rho^\gamma |v|^{p-2}\n vdsdx|\\
&\leq \e\int^t_0 \int_{\R^{N}} \rho|\n v|^2|v|^{p-2} dx ds+C(\e)\int^{t}_{0} \int_{\R^{N}}\rho^{2\gamma-1}|v|^{p-2} ds dx\\
&\leq \e\int^t_0 \int_{\R^{N}} \rho|\n v|^2|v|^{p-2} dx ds+C(\e)\int^{t}_{0} \int_{\R^{N}}\rho^{2\gamma-1}\rho^{\frac{p-2}{p}}|v|^{p-2} \rho^{-1+\frac{2}{p}}ds dx\\
&\leq \e\int^t_0 \int_{\R^{N}} \rho|\n v|^2|v|^{p-2} dx ds+C(\e)\int^{t}_{0} \int_{\R^{N}}(\rho^{2(\gamma-1+\frac{1}{p})}-1)\rho^{\frac{p-2}{p}}|v|^{p-2}ds dx\\
&\hspace{9cm}+C(\e)\int^{t}_{0} \int_{\R^{N}}\rho^{\frac{p-2}{p}}|v|^{p-2}ds dx\\
&\leq \e\int^t_0 \int_{\R^{N}} \rho|\n v|^2|v|^{p-2} dx ds+C(\e)\int^{t}_{0}( \int_{\R^{N}}\rho|v|^{p} dx)^{\frac{p-2}{p}} (\int_{\R^{N}}|\rho^{2(\gamma-1+\frac{1}{p})}-1|^{\frac{p}{2}} dx)^{\frac{2}{p}} ds\\
&\hspace{9cm}+C(\e)\int^{t}_{0} \int_{\R^{N}}\rho^{\frac{p-2}{p}}|v|^{p-2}ds dx\\
&\leq \e\int^t_0 \int_{\R^{N}} \rho|\n v|^2|v|^{p-2} dx ds+C(\e)\frac{\e_1(p-2)}{p}\int^{t}_{0} \int_{\R^{N}}\rho|v|^{p} dx ds+C(\e)\frac{2}{\e_1 p} \int^t_0\int_{\R^{N}}|\rho^{2(\gamma-1+\frac{1}{p})}-1|^{\frac{p}{2}} dx ds\\
&\hspace{9cm}+C(\e)\int^{t}_{0} \int_{\R^{N}}\rho^{\frac{p-2}{p}}|v|^{p-2}ds dx\\
\end{aligned}
\label{abtech}
\end{equation}
We have for $p\geq 4$ according to the H\"older's inequality that:
$$  \int_{\R^{N}}\rho|v|^{p-2}dx\leq\|\sqrt{\rho}v\|_{L^2}^{2(1-\theta)}\|\rho^{\frac{1}{p}}v\|_{L^p}^{p\theta},$$
with:
$$\theta=\frac{p-4}{p-2},\;1-\theta=\frac{2}{p-2}. $$
Then by H\"older's inequality we get:
$$  
\begin{aligned}
&\int^t_0 \int_{\R^{N}}\rho|v|^{p-2}dx\leq \|\sqrt{\rho}v \|_{L_t^\infty(L^2)}^{\frac{4}{p-2}} \int^t_0 \|\rho^{\frac{1}{p}}v(s)\|_{L^p}^{\frac{p(p-4)}{p-2}}ds.
\end{aligned}
$$
We deduce by Young inequality (and choosing $\e=\frac{1}{p}$) that:
\begin{equation}  
\begin{aligned}
&\int^t_0 \int_{\R^{N}}\rho|v|^{p-2}dx\leq \|\sqrt{\rho}v \|_{L_t^\infty(L^2)}^{\frac{4}{p-2}} \int^t_0 \|\rho^{\frac{1}{p}}v\|_{L^p}^{\frac{p(p-4)}{p-2}}ds,\\
&\leq \|\sqrt{\rho}v \|_{L_t^\infty(L^2)}^{\frac{4}{p-2}}  \int^t_0\big(\frac{\e(p-4)}{p-2} \|\rho^{\frac{1}{p}}v\|_{L^p}^{p}+\frac{2}{\e (p-2)}\big)ds,\\
&\leq \|\sqrt{\rho}v \|_{L_t^\infty(L^2)}^{\frac{4}{p-2}}  \big(\frac{(p-4)}{p-2} \int^t_0\frac{1}{p}\|\rho^{\frac{1}{p}}v(s,\cdot)\|_{L^p}^{p} ds+\frac{2p}{ (p-2)}t\big).
\end{aligned}
\label{etim1}
\end{equation}
We deduce that
\begin{equation}
\begin{aligned}
&|\int^{t}_{0} \int_{\R^{N}}|v|^{p-2}v\cdot\n a\rho^\gamma(s,x)dsdx|\\
&\leq \e\int^t_0 \int_{\R^{N}} \rho|\n v|^2|v|^{p-2} dx ds+C(\e)\frac{\e_1(p-2)}{p}\int^{t}_{0} \int_{\R^{N}}\rho|v|^{p} dx ds\\
&+C(\e)\frac{2}{\e_1 p} \int^t_0\int_{\R^{N}}|\rho^{2(\gamma-1+\frac{1}{p})}-1|^{\frac{p}{2}} dx ds\\&+C(\e) \|\sqrt{\rho}v \|_{L_t^\infty(L^2)}^{\frac{4}{p-2}}  \big(\frac{(p-4)}{p-2} \int^t_0\frac{1}{p}\|\rho^{\frac{1}{p}}v(s,\cdot)\|_{L^p}^{p} ds+\frac{2p}{ (p-2)}t\big).\\
\end{aligned}
\label{abtech1}
\end{equation}
We can conclude by a Gronwal argument for $p=5$ except we have to estimate when $p=5$:
$$\int^t_0\int_{\R^{N}}|\rho^{2(\gamma-1+\frac{1}{p})}-1|^{\frac{p}{2}} dx ds=\int^t_0\int_{\R^{N}}|\rho^{2(\gamma-\frac{4}{5})}-1|^{\frac{5}{2}} dx ds.$$
Now we recall that:
$$|\sqrt{\rho}-1|\leq|\rho-1|,$$
it implies that $\sqrt{\rho}-1$ belongs to $L^\infty_T(L^\gamma_2(\R^N))$. Furthermore if $\gamma\geq 2$ we deduce that  $\sqrt{\rho}-1$ is bounded in $L^\infty_T(L^2(\R^N))$ since $L^\gamma_2$ is embedded in this case in $L^2(\R^N)$. If $1\leq\gamma<2$, using the fact that $\n\sqrt{\rho}$ is bounded in $L^\infty_T(L^2(\R^N))$ we prove as in \cite{fH2} that $\sqrt{\rho}-1$ is in $L^\infty_T(L^2(\R^N))$. In particular it implies that $\sqrt{\rho}-1$ is bounded in $L^\infty_T(H^1(\R^N))$ for any $T>0$.\\
Now we can observe that
$$
\begin{aligned}
&|\rho^{2(\gamma-\frac{4}{5})}-1|\leq C|\rho-1|1_{\{|\rho-1|\leq\e\}}+C_1|\rho-1|^{2(\gamma-\frac{4}{5})}1_{\{|\rho-1|\geq\e\}},\\
&\leq C|\rho-1|1_{\{|\rho-1|\leq\e\}}+C_1(|\sqrt{\rho}-1|^2+2|\sqrt{\rho}-1|)^{2(\gamma-\frac{4}{5})}1_{\{|\rho-1|\geq\e\}}.
\end{aligned}$$
for $\e,C,C_1$ well chosen. We deduce by Sobolev embedding that:
$$\int^t_0\int_{\R^{N}}|\rho^{2(\gamma-\frac{4}{5})}-1|^{\frac{5}{2}} dx ds\leq t (1+\|\sqrt{\rho}-1\|_{L^\infty_t(H^1(\R^N))}^\alpha),$$
with $\alpha$ sufficiently large. We have finally prove that $\rho^{\frac{1}{5}}v$ belongs to $L^\infty_T(L^5(\R^N))$ for any $T>0$.
\\
\\
It implies that using the same bootstrap argument than in section \ref{section5}, we prove easily that $\rho$ belongs to $L^\infty_t(L^\infty(\R^N))$. We can now bound for any $p\geq 2$ the term (\ref{pression}) as follows::
$$|\int^{t}_{0} \int_{\R^{N}}|v|^{p-2}v\cdot\n a\rho^\gamma(s,x)dsdx|\leq C\|\rho^{\gamma-1}\|_{L^\infty_t(L^\infty)}\int^{t}_{0} \int_{\R^{N}}|{\rm div}(|v|^{p-2}v)|\rho(s,x)dsdx.$$
We can conclude the proof by using exactly the same argument than in section \ref{section4}.
\subsection{Case $N=3$}
We have now to deal with the term:
$$
\begin{aligned}
&|\int^{t}_{0} \int_{\R^{N}}|v|^{p-2}v\cdot\n a\rho^\gamma(s,x)dsdx|\leq C |\int^{t}_{0} \int_{\R^{N}}\rho^{\gamma-\frac{3}{2}+\frac{1}{p}}\rho^{\frac{p-1}{p}}|v|^{p-2}v\cdot\n\sqrt{\rho}(s,x)dsdx|
\end{aligned}
$$
\subsubsection{ Case $\frac{7}{6}<\gamma<\frac{7}{3}$}
Taking $p$ such that $p=3+\e$ with $\e>0$ small such that $\gamma-\frac{7+3\e}{2(3+\e)}>0$, by H\"older inequality and Sobolev embedding, we have:
$$
\begin{aligned}
&|\int^{t}_{0} \int_{\R^{N}}|v|^{p-2}v\cdot\n a\rho^\gamma(s,x)dsdx|\leq C |\int^{t}_{0} \int_{\R^{N}}\rho^{\gamma-\frac{7+3\e}{2(3+\e)}}\rho^{\frac{2+\e}{3+\e}}|v|^{1+\e}v\cdot\n\sqrt{\rho}(s,x)dsdx|\\
&\leq C\big(\| \n\sqrt{\rho}\|_{L^2_t(L^6(\R^N))}\|\rho^{\frac{1}{3+\e}}v\|^{2+\e}_{L^\infty_t(L^{3+\e}(\R^N))}
\|\rho^{\gamma-\frac{7+3\e}{2(3+\e)}}-1\|_{L^2_t(L^{\frac{6(3+\e)}{3-\e}})}\\
&\hspace{6cm}+\|\n\sqrt{\rho}\|_{L^\infty_t(L^2(\R^N))}\|\rho^{\frac{2+\e}{3+\e}}|v|^{2+\e}\|_{L^1_t(L^2(\R^N))}\big).
\end{aligned}
$$
Next we recall that via the energy inequality that we have $\sqrt{\rho}-1$ belongs to $L^\infty_t(H^1(\R^N))$, it yields by Sobolev embedding and energy inequality:
\begin{equation}
\n\rho^{\frac{\gamma}{2}}\in L^2_t(L^6(\R^N))\;\;\mbox{and}\;\;\rho-1\in L^\infty_t(L^3(\R^N)).
\label{energy}
\end{equation}
We have now to estimate:
$$\|\rho^{\gamma-\frac{7+3\e}{2(3+\e)}}-1\|_{L^2_t(L^{\frac{6(3+\e)}{3-\e}})}.$$
First we have for $\\alpha,C>0$:
$$|\rho^{\gamma-\frac{7+3\e}{2(3+\e)}}-1|1_{\{|\rho-1|\leq \alpha\}}\leq C|\rho-1|1_{\{|\rho-1|\leq \alpha\}}$$
In particular it implies that:
\begin{equation}
\|(\rho^{\gamma-\frac{7+3\e}{2(3+\e)}}-1)1_{\{|\rho-1|\leq \alpha\}}\|_{L^2_t(L^{\frac{6(3+\e)}{3-\e}})}\leq C(\e)\sqrt{t}\|\rho-1\|_{L^\infty_t(L^\gamma_2(\R^N))}.
\end{equation}
In a similar way, there exists $C_1>0$ such that:
$$|\rho^{\gamma-\frac{7+3\e}{2(3+\e)}}-1|1_{\{|\rho-1|\geq \alpha\}}\leq C_1
|\rho-1| ^{\gamma-\frac{7+3\e}{2(3+\e)}}     1_{\{|\rho-1|\geq \alpha\}}$$
We have then:
\begin{equation}
\|(\rho^{\gamma-\frac{7+3\e}{2(3+\e)}}-1)1_{\{|\rho-1|\geq \alpha\}}\|_{L^2_t(L^{\frac{6(3+\e)}{3-\e}})}\leq \|\rho-1\|_{L^{2(\gamma-\frac{7+3\e}{2(3+\e)})}_t(L^{(\gamma-\frac{7+3\e}{2(3+\e)})(\frac{6(3+\e)}{3-\e}})(\R^N))}.
\end{equation}
From (\ref{energy}), we deduce that if $\gamma>\frac{5}{3}$ (in order to ensure that $3\leq(\gamma-\frac{7+3\e}{2(3+\e)})(\frac{6(3+\e)}{3-\e})$ and $2(\gamma-\frac{7+3\e}{2(3+\e)})\geq 1$) and $\e$ sufficiently small:
\begin{equation}
\begin{aligned}
&(\rho-1)1_{\{|\rho-1|\geq \e\}}\in L_t^{\gamma}(L^{3\gamma}(\R^N))\cap L^\infty_t(L^{3}(\R^N)).
\end{aligned}
\label{interp11}
\end{equation}
By interpolation we have $(\rho-1)1_{\{|\rho-1|\geq \e\}}\in L^p_t(L^{(\gamma-\frac{7+3\e}{2(3+\e)})(\frac{6(3+\e)}{3-\e}})(\R^N))$:
$$
\begin{aligned}
&\frac{1}{6\gamma\frac{3+\e}{3-\e}-\frac{(7+3\e)3}{3-\e}}=\frac{\theta}{3\gamma}+\frac{1-\theta}{3}\\
&\frac{1}{p}=\frac{\theta}{\gamma}.
\end{aligned}
$$
We have then
$$
\begin{aligned}
&\theta=\frac{\gamma}{\gamma-1}\frac{2\gamma(3+\e)-(10+2\e)}{2\gamma(3+\e)-(7+3\e)}\;\;\mbox{and}\;\;p=(\gamma-1)\frac{2\gamma(3+\e)-(7+3\e)}{2\gamma(3+\e)-(10+2\e)}
\end{aligned}
$$
We observe now that:
$$p\geq2(\gamma-\frac{7+3\e}{2(3+\e)}).$$
It implies that we can bounded the term $\|\rho-1\|_{L^{2(\gamma-\frac{7+3\e}{2(3+\e)})}_t(L^{(\gamma-\frac{7+3\e}{2(3+\e)})(\frac{6(3+\e)}{3-\e}})(\R^N))}$. It conclude the proof when $\frac{5}{3}<\gamma<\frac{7}{3}$ since we have seen that $\rho^{\frac{1}{3+\e}}v\in L^\infty_T(L^{3+\e}(\R^N))$. It implies that using the same bootstrap argument than in section \ref{section5}, we prove easily that $\rho$ belongs to $L^\infty_t(L^\infty(\R^N))$. The rest of the proof follows the same line than in the section \ref{section4}.
\\
\\
Concerning the case $\frac{7}{6}<\gamma\leq\frac{5}{3}$, we have:
\begin{equation}
\|(\rho^{\gamma-\frac{7+3\e}{2(3+\e)}}-1)1_{\{|\rho-1|\geq \e\}}\|_{L^2_t(L^{\frac{6(3+\e)}{3-\e}})}\leq C(t+\|\rho-1\|_{L^{1}_t(L^{(\gamma-\frac{7+3\e}{2(3+\e)})(\frac{6(3+\e)}{3-\e}})(\R^N))}).
\end{equation}
We now use an interpolation argument as in (\ref{interp11}) except that we consider $(\rho-1)$ in $L^\infty_t(L^2(\R^N))$ (it is easy to observe that  $(\rho-1)$ is in $L^\infty_t(L^2(\R^N))$ by using the same arguments than \cite{fH2}). 
\subsubsection{ Case $1<\gamma<\frac{7}{6}$}
In this case we are going to estimate in an other way (\ref{pression}), more precisely we have for $C>0$:
\begin{equation}
\begin{aligned}
&|\int^{t}_{0} \int_{\R^{N}}|v|^{p-2}v\cdot\n a\rho^\gamma(s,x)dsdx|\leq C |\int^{t}_{0} \int_{\R^{N}}|v|^{p-2}|\n v|\rho^\gamma(s,x)dsdx|\\
&\leq\e\int^{t}_{0} \int_{\R^{N}}\rho|\n v|^2|v|^{p-2} dx ds+C(\e)\int^{t}_{0} \int_{\R^{N}}\rho^{2\gamma-1}|v|^{p-2}dx ds
\end{aligned}
\label{pression1}
\end{equation}
Because we have to deal with low frequencies, we are going to consider the case $p=4$. It gives using H\"older inequality:
\begin{equation}
\begin{aligned}
&|\int^{t}_{0} \int_{\R^{N}}|v|^{2}v\cdot\n a\rho^\gamma(s,x)dsdx|\leq\e\int^{t}_{0} \int_{\R^{N}}\rho|\n v|^2|v|^{2} dx ds+C(\e)\int^{t}_{0} \int_{\R^{N}}\rho^{2\gamma-2}\rho|v|^{2}dx ds\\
&\leq \e\int^{t}_{0} \int_{\R^{N}}\rho|\n v|^2|v|^{2} dx ds+C(\e)(\int^{t}_{0} \int_{\R^{N}}\rho|v|^4 dx ds+\frac{1}{2}\int^{t}_{0} \int_{\R^{N}}(\rho^{2\gamma-2}-1)^2 dx ds\\
&\hspace{5cm}+\frac{1}{2}\int^{t}_{0} \int_{\R^{N}}(\rho^{2\gamma-2}-1)^2(\sqrt{\rho}-1)^2 dx ds+t\|\sqrt{\rho}v\|_{L^\infty_t(L^2(\R^N))}.)
\end{aligned}
\label{pression1}
\end{equation}
We proceed as in the section \ref{section4} in order to prove that $\rho^{\frac{1}{4}}v$ belongs to $L^\infty_t(L^4(\R^N))$ via a bootstrap argument. It remains only to show that $\rho^{2\gamma-2}-1$ is bounded in $L^2_t(L^2(\R^N))$. We have for $\e,C,C_1$ suitably choose:
$$|\rho^{2\gamma-2}-1|\leq C|\rho-1| 1_{\{|\rho-1|\leq \e\}}+C_1|\rho-1|^{2\gamma-2}1_{\{|\rho-1|\geq \e\}}.$$
It suffices now to show that $(\rho-1)1_{\{|\rho-1|\geq \e\}}$ belongs to $L^{1}((0,t)\times\R^N)$ (since $4\gamma-4\leq1$). It is the case because $(\rho-1)1_{\{|\rho-1|\geq \e\}}$ is in $L^\infty_t(L^3(\Omega(t)))$ which is embedded in $L^{1}((0,t)\times\Omega(t))$ with $\Omega(t)=\{|\rho(t,\cdot)-1|\geq \e\}$ since $\Omega(t)$ is of finite measure.\\
It proves the fact that $\rho^{\frac{1}{4}}v$ is bounded in $L^\infty_T(L^4(\R^N))$ for any $T>0$. As in the section \ref{section5}, we show that $\rho$ is bounded in $L^\infty_T(L^\infty(\R^N))$ and we finish by proving that $\rho^{\frac{1}{p}}v$ is uniformly bounded in $p$ in $L^\infty_t(L^p(\R^N))$ for any $p\geq 2$. It concludes the proof of the theorem \ref{theo2} for general pressure.
\section{Appendix}
In this appendix, we only want to detail the computation on the Korteweg tensor.
\begin{lemme}
$${\rm div}K=\kappa{\rm div}(\rho\n\n\ln\rho)=\kappa{\rm div}(\rho D(\n\ln\rho)).$$
\end{lemme}
{\bf Proof:} By calculus, we obtain then:
\begin{equation}
\begin{aligned}
({\rm div}K)_{j}&=\big(\n\D\rho-{\rm div}(\frac{1}{\rho}\n\rho\otimes\n\rho)\big)_{j},\\
&=\p_{j}\D\rho-\frac{1}{\rho}\D\rho\,\p_{j}\rho-\frac{1}{2\rho}\p_{j}|\n\rho|^{2}+\frac{1}{\rho^{2}}|\n\rho|^{2}\p_{j}\rho,\\
\end{aligned}
\label{princip}
\end{equation}
Next we have:
$$\D\rho=\rho\D\ln\rho+\frac{1}{\rho}|\n\rho|^{2}.$$
We have then:
\begin{equation}
\begin{aligned}
\p_{j}\D\rho-\frac{1}{\rho}\D\rho\p_{j}\rho&=\p_{j}(\rho\,\D\ln\rho+\frac{1}{\rho}|\n\rho|^{2})
-\D\ln\rho\,\p_{j}\rho-\frac{1}{\rho^{2}}|\n\rho|^{2}\p_{j}\rho,\\
&=\rho\p_{j}\D\ln\rho+\frac{1}{\rho}\p_{j}(|\n\rho|^{2})-\frac{2}{\rho^{2}}|\n\rho|^{2}\p_{j}\rho.
\end{aligned}
\label{prin1}
\end{equation}
Putting the expression of (\ref{prin1}) in (\ref{princip}), we obtain:
\begin{equation}
\begin{aligned}
({\rm div}K)_{j}&=\p_{j}\D\rho+\frac{1}{2\rho}\p_{j}(|\n\rho|^{2})-\frac{1}{\rho^{2}}|\n\rho|^{2}\p_{j}\rho.
\end{aligned}
\label{principa}
\end{equation}
Next by calculus, we have:
\begin{equation}
\begin{aligned}
\frac{1}{2\rho}\p_{j}(|\n\rho|^{2})-\frac{1}{\rho^{2}}|\n\rho|^{2}\p_{j}\rho&=\sum_{i}(\p_{i}\ln\rho\p_{ij}\rho-(\p_{i}\ln\rho)^{2}\p_{j}\rho),\\
&=\sum_{i}\p_{i}\ln\rho\,\rho\p_{i,j}\ln\rho,\\
&=\frac{\rho}{2}\n(|\ln\rho|^{2})_{j}.
\end{aligned}
\label{imp}
\end{equation}
Finally by using (\ref{imp}) and (\ref{principa}), we obtain:
$$
\begin{aligned}
{\rm div}K=\rho(\n\D(\ln\rho)+\frac{\rho}{2}\n(|\n\ln\rho|^{2})).
\label{cap}
\end{aligned}
$$
We now want to prove that we can rewrite (\ref{cap}) under the form of a viscosity tensor. To see this, we have:
$$
\begin{aligned}
{\rm div}(\rho\n(\n \ln\rho))_{j}&=\sum_{i}\p_{i}(\rho\p_{ij}\ln\rho),\\
&=\sum_{i}[\p_{i}\rho\p_{ij}\ln\rho+\rho\p_{iij}\ln\rho],\\
&=\rho(\D\n \ln\rho)_{j}+\sum_{i}\rho\p_{i}\ln\rho\p_{j}\p_{i}\ln\rho),\\
&=\rho(\D\n \ln\rho)_{j}+\frac{\rho}{2}(\n(|\n\ln\rho|^{2}))_{j},\\
&={\rm div}K.
\end{aligned}
$$
We have then:
$${\rm div}K=\kappa{\rm div}(\rho\n\n\ln\rho)=\kappa{\rm div}(\rho D(\n\ln\rho)).$$
$\blacksquare$

\end{document}